\documentclass[preprint,fleqn,10pt]{elsarticle}
\usepackage[T1]{fontenc}
\usepackage[utf8]{inputenc}
\usepackage{mathptmx}
\usepackage[scaled=0.94]{helvet}
\usepackage[paperwidth=155mm,paperheight=235mm,left=16.65mm,right=16.1mm,
  top=22mm,bottom=18.5mm,headheight=11pt,headsep=12pt]{geometry}
\usepackage{amssymb,mathtools,amsthm,microtype}
\usepackage[colorlinks=true,linkcolor=blue,citecolor=blue,urlcolor=blue]{hyperref}
\hypersetup{pdftitle={Approximate Lr-factorization of lp-strictly singular operators on Lp},
  pdfauthor={Gang Cai, Jinghao Huang, Yin Zhou}}
\input{glyphtounicode}
\allowdisplaybreaks[1]
\numberwithin{equation}{section}
\newtheorem{theorem}{Theorem}[section]
\newtheorem{lemma}{Lemma}[section]
\newtheorem{proposition}{Proposition}[section]
\newtheorem{corollary}{Corollary}[section]
\theoremstyle{definition}
\newtheorem{problem}{Problem}[section]
\DeclareMathOperator{\dist}{dist}

\begin{document}
\begin{frontmatter}
\title{Approximate \(L_r\)-factorization
  of \(\ell_p\)-strictly singular operators on \(L_p\)\tnoteref{funding}}

\author[cqnu]{Gang Cai\corref{cor}}
\ead{caigang-aaaa@163.com}
\author[hit]{Jinghao Huang}
\ead{jinghao.huang@hit.edu.cn}
\author[cqnu]{Yin Zhou}
\ead{zhouyin0330@163.com}
\address[cqnu]{School of Mathematical Sciences, Chongqing Normal University,
  Chongqing 401331, China}
\address[hit]{Institute for Advanced Study in Mathematics of HIT,
  Harbin Institute of Technology, Harbin 150001, China}
\cortext[cor]{Corresponding author.}
\tnotetext[funding]{Gang Cai is supported by the National Natural Science
  Foundation of China (grant no. 12571141) and the Chongqing Talent Support
  Program (grant no. cstc2024ycjh-bgzxm0121).
  Jinghao Huang was supported by the National Natural
  Science Foundation of China (grant nos. 12301160, 12471134 and 12671159).}

\begin{abstract}
For \(1<p<2\), we study norm approximation of operators on \(L_p(0,1)\)
by operators that factor through \(L_r(0,1)\), where \(p<r<2\).
Using a sharp factorization estimate for finite diagonal blocks, we show
that the corresponding norm-closed factorization ideals form a strictly
decreasing chain as \(r\) increases. We also construct a strictly singular
operator \(T_p\) at positive distance from the closed linear span generated
by these ideals and all finitely strictly singular operators on \(L_p(0,1)\).
In particular, \(T_p\) is \(\ell_p\)-strictly singular, yet its distance from
the \(L_r\)-factorization ideals has a positive lower bound independent of
\(r\). This answers a question due to Johnson and Schechtman in the negative.
\end{abstract}

\end{frontmatter}
\thispagestyle{plain}

\section{Introduction}
We study norm approximation of operators on \(L_p(0,1)\) by operators
factoring through \(L_r(0,1)\), where \(1<p<r<2\). Write \(L_s=L_s(0,1)\)
with Lebesgue measure. Let \(E\) and \(F\) be Banach spaces. A bounded
linear operator \(T:E\to F\) is \emph{bounded below} on a closed subspace
\(X\subseteq E\) if \(\|Tx\|\ge c\|x\|\) for all \(x\in X\) and some
\(c>0\), or equivalently, if \(T|_X\) is an isomorphism onto its range.
For a Banach space \(Z\), the operator \(T\) is \emph{\(Z\)-strictly
singular} if it is bounded below on no closed subspace isomorphic to
\(Z\) \cite{JS}. An operator is \emph{strictly singular} if it is bounded
below on no infinite-dimensional closed subspace \cite{Kato,MT}.

An operator \(S:E\to F\) \emph{factors through \(L_r\)} if \(S=BA\)
for bounded operators \(A:E\to L_r\) and \(B:L_r\to F\). Write
\(\mathcal G_r(E,F)\) for this class and \(\overline{\mathcal G_r(E,F)}\)
for its closure in the operator norm. For a nonempty class \(\mathcal C\)
of operators from \(E\) to \(F\), put
\(\dist(T,\mathcal C)=\inf_{S\in\mathcal C}\|T-S\|\).

Answering a question of Saksman and Tylli \cite{ST}, Johnson and
Schechtman \cite[Theorem~1]{JS} characterized the weakly compact
multiplication operators \(S\mapsto ASB\) on \(\mathcal L(L_p)\).
They also proved that if an \(\ell_p\)-strictly singular operator \(T\) on
\(L_p\), \(1<p<2\), is bounded below on a closed subspace \(X\), then
\(X\) embeds isomorphically into \(L_r\) for every \(p<r<2\)
\cite[Theorem~4]{JS}. Such an embedding is necessary for approximation
of \(T\) through \(L_r\). Indeed, if \(T|_X\) has lower bound \(c>0\)
and \(S=BA\) factors through \(L_r\) with \(\|T-S\|<c\), then
\[
 (c-\|T-S\|)\|x\|\le\|BAx\|\le\|B\|\,\|Ax\| \qquad(x\in X),
\]
so the first factor \(A\) embeds \(X\) into \(L_r\).
Johnson and Schechtman asked whether \(\ell_p\)-strict singularity
also suffices for approximation on the whole of \(L_p\).

\begin{problem}[{Johnson and Schechtman \cite[Problem 2]{JS}}]\label{prob:1}
Given $1<p<r<2$ and an $\ell_{p}$-strictly singular operator $T: L_{p} \rightarrow L_{p}$, must $T$ belong to $\overline{\mathcal{G}_{r}\left(L_{p}, L_{p}\right)}$?
\end{problem}

The following theorem gives a negative answer with a single operator
for all intermediate exponents.

\begin{theorem}\label{res:1}
For each $1<p<2$, there exist an $\ell_{p}$-strictly singular operator $T_{p}: L_{p} \rightarrow L_{p}$ and a constant $\eta_{p}>0$ such that
\begin{equation}
\inf _{p<r<2} \operatorname{dist}\left(T_{p}, \overline{\mathcal{G}_{r}\left(L_{p}, L_{p}\right)}\right) \geq \eta_{p} . \label{eq:1}
\end{equation}
\end{theorem}

The \(\ell_p\)-strict singularity condition in Problem~\ref{prob:1}
is necessary. Indeed, if an operator factoring through \(L_r\) were bounded below on a copy
of \(\ell_p\), its first factor would induce an embedding
\(V:\ell_p\to L_r\). For \(p<r\le2\), a lower bound \(b>0\) for \(V\)
and the type-\(r\) inequality in \(L_r\), applied to \(Ve_1,\ldots,Ve_n\)
with independent signs \(\varepsilon_j\), would give
\[
 b n^{1/p}\le
 \left(\mathbb E\Big\|\sum_{j=1}^n\varepsilon_j V e_j\Big\|_r^r\right)^{1/r}
 \le K_r\Big(\sum_{j=1}^n\|V e_j\|_r^r\Big)^{1/r}
 \le K_r\|V\|n^{1/r}.
\]
Since \(1/p>1/r\), this is impossible for large \(n\). Boundedness
below on a fixed subspace persists under small norm perturbations, so
every norm limit of such operators is also \(\ell_p\)-strictly singular.

At the endpoint \(r=2\), Johnson and Schechtman \cite[Section~4]{JS}
disproved Tylli's conjecture, even in its weaker form
allowing an isometric embedding of the range into \(\ell_\infty\).
For biased-coin convolution operators on the Cantor group,
Pisier \cite[Theorem~1]{Pisier} established norm approximation by
operators factoring through a Hilbert space. The endpoint counterexample
leaves open the possibility of approximation through an intermediate
\(L_r\), since for \(p<r<2\) the closed
\(L_2\)-factorization ideal is contained in the closed
\(L_r\)-factorization ideal (see the proof of Corollary~\ref{cor:jumps}).

Let \(q=p'\) and \(s=r'\), so that \(2<s<q\). Taking adjoints turns an
\(L_r\)-factorization into an \(L_s\)-factorization and preserves the approximation
error. We use the weighted diagonal model of Johnson and Schechtman
\cite[Section~4]{JS}, in which decay of the diagonal coefficients to zero
ensures \(\ell_q\)-strict singularity. Proposition~\ref{res:2} determines
the minimum product of the factor norms required to approximate a finite
block through \(L_s\) within a prescribed error. We choose a single
sequence of block weights so that this minimum tends to infinity along
the blocks for every \(2<s<q\) and every fixed error tolerance
\(0\leq\varepsilon<1\). In Section~\ref{sec:main} we realize the resulting
diagonal \(D\) as a compression of \(T_p^*\). The distance estimate for \(D\)
in Lemma~\ref{res:4} and the compression bounds in Proposition~\ref{res:5}
then give the same positive lower bound for the distance from \(T_p\)
to each of the closed \(L_r\)-factorization ideals, \(p<r<2\).

Varying the weights also allows us
to prescribe either \((p,a)\) or \((p,a]\) as the set of intermediate
exponents for which an operator factors, for any \(p<a<2\)
(Theorem~\ref{res:7}). It follows that
the closed ideal at \(a\) strictly contains the closure of the union
of the ideals at larger exponents up to \(2\), and is strictly contained
in the intersection of the ideals at exponents between \(p\) and \(a\)
(Corollary~\ref{cor:jumps}). For the broader structure of closed
ideals in \(\mathcal L(L_p)\), see \cite{JSApprox,SchechtmanICM}.

The operator constructed for Theorem~\ref{res:1} is bounded below on a
copy of \(\ell_2\). This raises the question of whether a strictly
singular operator can also give a counterexample. An operator on \(L_p\) is strictly singular
if and only if it is both \(\ell_p\)-strictly singular and
\(\ell_2\)-strictly singular \cite{KP,Weis}. Strictly singular operators
can be noncompact \cite{GT}. For \(p>2\), \(\ell_2\)-strict
singularity on \(L_p\) implies factorization through \(\ell_p\)
\cite{J,KP}. In our range \(1<p<2\), Mathieu and Tradacete
\cite[Theorem~6.2]{MT} proved that every strictly singular operator on
\(L_p\) factors through \(\ell_p\), whereas Johnson and Schechtman
\cite[Remark~4.3]{JSIdeals} showed that some such operators cannot be
approximated through \(L_2\).

Section~\ref{sec:strict} constructs a strictly singular operator at
positive distance from the closed linear span of all intermediate
factorizations and all finitely strictly singular operators
(Theorems~\ref{res:8} and~\ref{res:9}). The quotient of the strictly
singular operators by those in this closed span contains a complemented
copy of \(\ell_\infty/c_0\) (Corollary~\ref{cor:quotient}).
The second counterexample is bounded below, with a common positive
constant, on subspaces of arbitrarily large finite dimension
by \eqref{eq:10}. Finite strict singularity rules out exactly this
behavior. Section~\ref{sec:finite} shows that finite strict singularity
does imply approximation for the diagonal model of Theorem~\ref{res:10}
and for the operators with martingale-difference columns in
Theorem~\ref{res:11}. Problem~\ref{prob:2} asks whether finite strict
singularity suffices for every operator on \(L_p\).

\section{Proof of the main theorem}\label{sec:main}

\subsection{A sharp factorization estimate for finite diagonal blocks}
Let \(q\) and \(s\) be real exponents satisfying \(2<s<q<\infty\).
For the application to Theorem~\ref{res:1}, fix \(p\in(1,2)\) and set
\(q=p'=p/(p-1)\). As \(r\) ranges over \((p,2)\), its conjugate exponent
\(s=r'=r/(r-1)\) ranges over the entire interval \((2,q)\).
If an approximant factors through \(L_r\), its adjoint factors through \(L_s\).
For a scalar diagonal block, the normalized trace yields a lower bound on the product of the factor norms for
every approximant. An explicit factorization attains this bound. Schlumprecht and Zs\'ak \cite{SZ} also use averages of matrix coefficients to separate operator ideals.

For a bounded operator $S: E \rightarrow F$ define
\[
\gamma_{s}(S)=\inf \left\{\|A\|\|B\|: S=B A, A: E \rightarrow L_{s}, B: L_{s} \rightarrow F\right\},
\]
with $\gamma_{s}(S)=+\infty$ if there is no such factorization.
We have $\|S\| \leq \gamma_{s}(S)$. For contractions $J$ and $P$, the identity $P S J=(P B)(A
J)$ gives $\gamma_{s}(P S J) \leq \gamma_{s}(S)$. In particular, a global factorization gives
a common upper bound on the factorization norms of all coordinate compressions.

\begin{proposition}\label{res:2}
Let $2<s<q<\infty$, let $N$ be a positive integer, and set
\[
w=N^{1 / q-1 / 2}, \quad w \leq \delta \leq 1, \quad v=w / \delta .
\]

Let $E=\ell_{q}^{N}$, and let $Y$ be the same coordinate space with norm
\[
\|y\|_{Y}=\max \left\{\|y\|_{q}, v\|y\|_{2}\right\} .
\]

Put $D_{N}=\delta I: E \rightarrow Y$. For every $0 \leq \varepsilon<1$,
\begin{equation}
\inf _{\substack{S: E \rightarrow Y \\\left\|D_{N}-S\right\| \leq \varepsilon}} \gamma_{s}(S)=(1-\varepsilon) \max \left\{1, \delta N^{1 / s-1 / q}\right\} . \label{eq:2}
\end{equation}

The infimum is attained by $S=(1-\varepsilon) D_{N}$.
\end{proposition}
\begin{proof}
For $x \in E$, H\"older's inequality gives
\[
\|x\|_{2}^{2} \leq\left(\sum_{j=1}^{N}\left|x_{j}\right|^{q}\right)^{2 / q} N^{1-2 / q}, \quad \frac{\|(1, \ldots, 1)\|_{2}}{\|(1, \ldots, 1)\|_{q}}=N^{1 / 2-1 / q} .
\]

Consequently,
\[
\begin{aligned}
\|I: E \rightarrow Y\| & =\sup _{\|x\|_{q}=1} \max \left\{1, v\|x\|_{2}\right\}=\max \left\{1, v N^{1 / 2-1 / q}\right\}=\delta^{-1}, \\
\left\|D_{N}\right\| & =\delta\|I: E \rightarrow Y\|=1 .
\end{aligned}
\]

For $S: E \rightarrow Y$, let $\tau(S)=N^{-1} \operatorname{tr} S$, where $\operatorname{tr}
S=\sum_{j=1}^{N}\left(S e_{j}\right)_{j}$ is the trace of the matrix of $S$ in the common
coordinate basis of $E$ and $Y$. Let $D_{\sigma}=\operatorname{diag}\left(\sigma_{1}, \ldots,
\sigma_{N}\right)$ for $\sigma \in\{-1,1\}^{N}$. The average over all sign changes satisfies
\[
\left[2^{-N} \sum_{\sigma} D_{\sigma}^{-1} S D_{\sigma}\right]_{i j}=S_{i j} 2^{-N} \sum_{\sigma} \sigma_{i} \sigma_{j}= \begin{cases}S_{i i}, & i=j, \\ 0, & i \neq j .\end{cases}
\]

For a permutation $\pi$, write $P_{\pi} e_{j}=e_{\pi(j)}$. For each coordinate vector $e_{j}$,
\[
P_{\pi}^{-1} \operatorname{diag}\left(S_{11}, \ldots, S_{N N}\right) P_{\pi} e_{j}=S_{\pi(j), \pi(j)} e_{j} .
\]

For fixed $j$ and $k$, there are $(N-1)!$ permutations with $\pi(j)=k$. Hence the averaged $j$-th diagonal entry is
\[
\frac{1}{N!} \sum_{\pi \in \mathfrak{S}_{N}} S_{\pi(j), \pi(j)}=\frac{(N-1)!}{N!} \sum_{k=1}^{N} S_{k k}=\tau(S) .
\]

Thus $|G|^{-1} \sum_{U \in G} U^{-1} S U=\tau(S) I$ for the group $G=\left\{D_{\sigma}
P_{\pi}\right\}$ of signed permutations. Each $U \in G$ acts isometrically on $E$ and on $Y$.

Since $U^{-1} D_{N} U=D_{N}$ for every $U \in G$,
\[
\begin{aligned}
\delta^{-1}|\delta-\tau(S)| & =\left\|D_{N}-\tau(S) I\right\| \\
& =\left\|\frac{1}{|G|} \sum_{U \in G} U^{-1}\left(D_{N}-S\right) U\right\| \\
& \leq \frac{1}{|G|} \sum_{U \in G}\left\|U^{-1}\left(D_{N}-S\right) U\right\|=\left\|D_{N}-S\right\| .
\end{aligned}
\]

For $\left\|D_{N}-S\right\| \leq \varepsilon$, the reverse triangle inequality gives
\begin{equation}
|\tau(S)| \geq \delta-|\delta-\tau(S)| \geq \delta-\varepsilon \delta=(1-\varepsilon) \delta . \label{eq:3}
\end{equation}

Fix a factorization $S=B A$ through $L_{s}$ and put $f_{j}=A e_{j}$. Since $\left|y_{j}\right| \leq\|y\|_{q} \leq\|y\|_{Y}$, each diagonal entry satisfies
\[
\left|S_{j j}\right|=\left|\left(S e_{j}\right)_{j}\right|=\left|\left(B f_{j}\right)_{j}\right| \leq\left\|B f_{j}\right\|_{Y} \leq\|B\|\left\|f_{j}\right\|_{s} .
\]

With $s^{\prime}=s /(s-1)$, H\"older's inequality yields
\begin{align}
|\tau(S)| & \leq \frac{\|B\|}{N} \sum_{j=1}^{N}\left\|f_{j}\right\|_{s} \notag\\
& \leq \frac{\|B\|}{N}\left(\sum_{j=1}^{N}\left\|f_{j}\right\|_{s}^{s}\right)^{1 / s}\left(\sum_{j=1}^{N} 1^{s^{\prime}}\right)^{1 / s^{\prime}}=\|B\| N^{-1 / s}\left(\sum_{j=1}^{N}\left\|f_{j}\right\|_{s}^{s}\right)^{1 / s} . \label{eq:4}
\end{align}

Let $\sigma$ be uniform on $\{-1,1\}^{N}$. For almost every $t$ put $a_{j}=f_{j}(t)$ and $X_{\sigma}(t)=\sum_{j=1}^{N} \sigma_{j} a_{j}$. The sign average above gives
\[
\mathbb{E}_{\sigma}\left|X_{\sigma}(t)\right|^{2}=\sum_{i, j=1}^{N} a_{i} \overline{a_{j}} \mathbb{E}_{\sigma}\left(\sigma_{i} \sigma_{j}\right)=\sum_{j=1}^{N}\left|a_{j}\right|^{2} .
\]

Writing $M=\left(\sum_{j}\left|a_{j}\right|^{2}\right)^{1 / 2}$ and applying Jensen's inequality to $u \mapsto u^{s / 2}$, we obtain
\[
\sum_{j=1}^{N}\left|a_{j}\right|^{s} \leq M^{s-2} \sum_{j=1}^{N}\left|a_{j}\right|^{2}=\left(\mathbb{E}_{\sigma}\left|X_{\sigma}(t)\right|^{2}\right)^{s / 2} \leq \mathbb{E}_{\sigma}\left|X_{\sigma}(t)\right|^{s} .
\]

Since $A \sigma=\sum_{j} \sigma_{j} f_{j}$ and $\|\sigma\|_{q}=N^{1 / q}$,
\[
\begin{aligned}
\sum_{j=1}^{N}\left\|f_{j}\right\|_{s}^{s} & \leq \int_{0}^{1} 2^{-N} \sum_{\sigma}\left|X_{\sigma}(t)\right|^{s} d t=2^{-N} \sum_{\sigma} \int_{0}^{1}\left|X_{\sigma}(t)\right|^{s} d t \\
& =\mathbb{E}_{\sigma}\|A \sigma\|_{s}^{s} \leq\|A\|^{s} \mathbb{E}_{\sigma}\|\sigma\|_{q}^{s}=\|A\|^{s} N^{s / q}
\end{aligned}
\]

Substitution in \eqref{eq:4} gives
\begin{equation}
|\tau(S)| \leq\|B\| N^{-1 / s}\left(\|A\|^{s} N^{s / q}\right)^{1 / s}=\|A\|\|B\| N^{1 / q-1 / s} . \label{eq:5}
\end{equation}

If $\left\|D_{N}-S\right\| \leq \varepsilon$, taking the infimum over factorizations in \eqref{eq:5} and using \eqref{eq:3} gives
\[
\begin{aligned}
& \gamma_{s}(S) \geq N^{1 / s-1 / q}|\tau(S)| \geq(1-\varepsilon) \delta N^{1 / s-1 / q}, \\
& \gamma_{s}(S) \geq\|S\| \geq\left\|D_{N}\right\|-\left\|D_{N}-S\right\| \geq 1-\varepsilon .
\end{aligned}
\]

Their maximum is the lower bound in \eqref{eq:2}.
For the upper bound, factor $D_{N}$ as
\[
\ell_{q}^{N} \xrightarrow{I} \ell_{s}^{N} \xrightarrow{\delta I} Y .
\]

Denote these factors by $A_{0}$ and $B_{0}$, respectively. Then
\[
\begin{aligned}
\left\|A_{0}\right\| & =\sup _{\|x\|_{q}=1}\|x\|_{s}=N^{1 / s-1 / q}, \\
\left\|B_{0}\right\| & =\delta \max \left\{\sup _{\|y\|_{s}=1}\|y\|_{q}, v \sup _{\|y\|_{s}=1}\|y\|_{2}\right\}=\delta \max \left\{1, v N^{1 / 2-1 / s}\right\} \\
& =\max \left\{\delta, N^{1 / q-1 / s}\right\}, \\
\left\|B_{0}\right\|\left\|A_{0}\right\| & =\max \left\{1, \delta N^{1 / s-1 / q}\right\} .
\end{aligned}
\]

In the second line, the first supremum is attained at $e_{1}$ and the second at $N^{-1 / s}(1, \ldots, 1)$.
Partition $(0,1)$ into $N$ intervals $J_{j}$ of equal measure and set
\[
J x=N^{1 / s} \sum_{j=1}^{N} x_{j} \mathbf{1}_{J_{j}}, \quad(Q f)_{j}=N^{1 / s^{\prime}} \int_{J_{j}} f(t) d t
\]

Then $(Q J x)_{j}=N\left|J_{j}\right| x_{j}=x_{j}$, and
\[
\|J x\|_{s}^{s}=N \sum_{j=1}^{N}\left|J_{j}\right|\left|x_{j}\right|^{s}=\|x\|_{s}^{s}, \quad\left|(Q f)_{j}\right| \leq\left(N\left|J_{j}\right|\right)^{1 / s^{\prime}}\left(\int_{J_{j}}|f(t)|^{s} d t\right)^{1 / s} .
\]

Raising both sides of the estimate for \(\lvert(Qf)_j\rvert\) to the \(s\)th power and summing over \(j\) gives
\(\|Qf\|_{\ell_s^N}\leq\|f\|_{L_s}\). Thus
$D_{N}=\left(B_{0} Q\right)\left(J A_{0}\right)$ factors through $L_{s}$. For
$S_{\varepsilon}=(1-\varepsilon) D_{N}$ this gives
\[
\begin{aligned}
\left\|D_{N}-S_{\varepsilon}\right\| & =\varepsilon\left\|D_{N}\right\|=\varepsilon, \\
\gamma_{s}\left(S_{\varepsilon}\right) & \leq(1-\varepsilon)\left\|B_{0} Q\right\|\left\|J A_{0}\right\| \leq(1-\varepsilon) \max \left\{1, \delta N^{1 / s-1 / q}\right\} .
\end{aligned}
\]

The lower bound already proved gives equality in \eqref{eq:2}.
\end{proof}

\subsection{Separation from intermediate \texorpdfstring{$L_{s}$}{Ls}-factorizations}
We use the weighted diagonal of Johnson and Schechtman \cite[Section 4]{JS}, who observed that
the condition $w_{j} / v_{j} \rightarrow 0$ implies $\ell_{q}$-strict singularity.
We choose the block lengths
$N_{k}$ and weight ratios $\delta_{k}$ so that $\delta_{k} \rightarrow 0$ while $\delta_{k}
N_{k}^{1 / s-1 / q} \rightarrow \infty$ for every fixed $2<s<q$. For any fixed
error less than one, Proposition~\ref{res:2} then gives a lower bound on
the block factorization norms that tends to infinity with $k$. This yields
the distance estimate in Lemma~\ref{res:4}.

Fix $2<q<\infty$. For weights $u=\left(u_{j}\right)$ with $0<u_{j} \leq 1$, let $X_{q, u}$ be the completion of the finitely supported sequences in the norm
\[
\|x\|_{q, u}=\max \left\{\|x\|_{q},\|u x\|_{2}\right\}, \quad(u x)_{j}=u_{j} x_{j} .
\]

Partition $\mathbb{N}$ into consecutive finite blocks $I_{k}$, indexed by $k \geq k_{0}$, and take
\begin{equation}
\left|I_{k}\right|=N_{k}=2^{k}, \quad \delta_{k}=\frac{1}{k}, \quad a_{k}=N_{k}^{1 / q-1 / 2}, \quad b_{k}=\frac{a_{k}}{\delta_{k}} . \label{eq:6}
\end{equation}

Since $q>2$, we can choose $k_{0}$, depending only on $q$, so that $a_{k} \leq \delta_{k}$ for
every $k \geq k_{0}$. Define $w_{j}=a_{k}, v_{j}=b_{k}$, and $d_{j}=\delta_{k}$ for $j \in
I_{k}$. Then $0<w_{j} \leq v_{j} \leq 1$. Define the diagonal map
\[
D: X_{q, w} \longrightarrow X_{q, v}, \quad(D x)_{j}=d_{j} x_{j} .
\]

The identity $v_{j} d_{j}=w_{j}$ gives $\|D x\|_{q, v}=\max \left\{\|d x\|_{q},\|w
x\|_{2}\right\} \leq\|x\|_{q, w}$. On $I_{k}$, the domain norm is exactly the
$\ell_{q}^{N_{k}}$ norm, since $a_{k}\|x\|_{2} \leq\|x\|_{q}$ on the block. The range norm is
$\max \left\{\|y\|_{q}, b_{k}\|y\|_{2}\right\}$. Hence the block restriction of $D$ is
$D_{N_{k}}$ from Proposition~\ref{res:2}. In particular, the vector $x=N_{k}^{-1 / q}
\mathbf{1}_{I_{k}}$ has $\|x\|_{q}=\|w x\|_{2}=1$. Thus every block has a unit vector whose
norm $D$ preserves, even though $\delta_{k} \rightarrow 0$, and $\|D\|=1$.

\begin{lemma}\label{res:3}
The operator $D$ is $\ell_{q}$-strictly singular.
\end{lemma}
\begin{proof}
Suppose $D$ is bounded below by $c>0$ on a closed subspace $X \subseteq X_{q, w}$ isomorphic
to $\ell_{q}$. Let $\mathbb{K}$ be the scalar field. Choose $m$ so that $d_{j} \leq c / 2$ for
$j>m$, and set
\[
R: X \rightarrow \mathbb{K}^{m}, \quad R x=\left(x_{1}, \ldots, x_{m}\right), \quad X_{0}=\operatorname{ker} R, \quad F=R(X) .
\]

For $z \in X_{0}$,
\[
\|d z\|_{q}=\left(\sum_{j>m}\left|d_{j} z_{j}\right|^{q}\right)^{1 / q} \leq(c / 2)\|z\|_{q} \leq(c / 2)\|z\|_{q, w} .
\]

Since $c\|z\|_{q, w} \leq\|D z\|_{q, v}=\max \left\{\|d z\|_{q},\|w z\|_{2}\right\}$, the second term gives the lower bound whenever $z \neq 0$. Thus
\[
W: X_{0} \rightarrow \ell_{2}, \quad W z=w z, \quad c\|z\|_{q, w} \leq\|W z\|_{2} \leq\|z\|_{q, w} .
\]

In particular $H=W\left(X_{0}\right)$ is closed and $\left\|W^{-1}: H \rightarrow X_{0}\right\| \leq 1 / c$.
Equip $F$ with its Euclidean norm, so that $\|R x\|_{2} \leq \sqrt{m}\|x\|_{q} \leq
\sqrt{m}\|x\|_{q, w}$. Choose a basis $y_{1}, \ldots, y_{d}$ of $F$ and $u_{i} \in X$ with $R
u_{i}=y_{i}$. Define
\[
S: F \rightarrow X, \quad S\left(\sum_{i=1}^{d} a_{i} y_{i}\right)=\sum_{i=1}^{d} a_{i} u_{i} .
\]

Then $R S=I_{F}$, and $S$ is bounded since $F$ is finite-dimensional. For $x \in X$, the identity $R(x-S R x)=0$ gives the decomposition
\[
x=\underbrace{x-S R x}_{\in X_{0}}+\underbrace{S R x}_{\in S(F)} .
\]

An isomorphism $U: X \rightarrow H \oplus_{2} F$ and its inverse are given by
\[
U x=(W(x-S R x), R x), \quad U^{-1}(h, y)=W^{-1} h+S y .
\]

Choose an isomorphism $J: \ell_{q} \rightarrow X$ and put $\mathcal{H}=H \oplus_{2} F$. The isomorphism $V=U J$ : $\ell_{q} \rightarrow \mathcal{H}$ satisfies
\[
b\|x\|_{q} \leq\|V x\|_{\mathcal{H}} \leq\|V\|\|x\|_{q} \quad\left(x \in \ell_{q}\right), \quad b=\left\|V^{-1}\right\|^{-1}>0 .
\]

Fix $n \geq 1$ and write $h_{j}=V e_{j}$ for $1 \leq j \leq n$. Since each $e_{j}$ has norm
one, $\sum_{j=1}^{n}\left\|h_{j}\right\|_{\mathcal{H}}^{2} \geq n b^{2}$. To estimate the same
sum from above, let $\sigma_{1}, \ldots, \sigma_{n}$ be independent uniform signs. Expanding
the squared Hilbert norm gives
\[
\mathbb{E}_{\sigma}\left\|\sum_{j=1}^{n} \sigma_{j} h_{j}\right\|_{\mathcal{H}}^{2}=\sum_{i, j=1}^{n} \mathbb{E}_{\sigma}\left(\sigma_{i} \sigma_{j}\right)\left\langle h_{i}, h_{j}\right\rangle_{\mathcal{H}}=\sum_{j=1}^{n}\left\|h_{j}\right\|_{\mathcal{H}}^{2} .
\]

Indeed, $\sigma_{j}^{2}=1$, whereas for $i \neq j$ independence gives
$\mathbb{E}_{\sigma}\left(\sigma_{i} \sigma_{j}\right)=\left(\mathbb{E}_{\sigma}
\sigma_{i}\right)\left(\mathbb{E}_{\sigma} \sigma_{j}\right)=0$. For every choice of signs,
\[
\left\|\sum_{j=1}^{n} \sigma_{j} h_{j}\right\|_{\mathcal{H}}=\left\|V \sum_{j=1}^{n} \sigma_{j} e_{j}\right\|_{\mathcal{H}} \leq\|V\|\left(\sum_{j=1}^{n}\left|\sigma_{j}\right|^{q}\right)^{1 / q}=\|V\| n^{1 / q} .
\]

Squaring and averaging yields
\[
n b^{2} \leq \sum_{j=1}^{n}\left\|h_{j}\right\|_{\mathcal{H}}^{2} \leq\|V\|^{2} n^{2 / q}, \quad n^{1-2 / q} \leq(\|V\| / b)^{2} .
\]

The last inequality must hold for every $n$ with the same $V$ and $b$. Since $1-2 / q=(q-2) /
q>0$, any integer $n>(\|V\| / b)^{2 q /(q-2)}$ violates it. Thus $D$ cannot be bounded below
on a copy of $\ell_{q}$.
\end{proof}

\begin{lemma}\label{res:4}
For every $2<s<q$,
\begin{equation}
\operatorname{dist}\left(D, \overline{\mathcal{G}_{s}\left(X_{q, w}, X_{q, v}\right)}\right)=1 . \label{eq:7}
\end{equation}

\end{lemma}
\begin{proof}
Fix $2<s<q$. Suppose that an operator $R=B A$, with $A: X_{q, w} \rightarrow L_{s}$ and $B:
L_{s} \rightarrow X_{q, v}$, satisfies $\varepsilon=\|D-R\|<1$. Put $M=\|A\|\|B\|$.
For each \(k\ge k_0\), index the \(N_k\) coordinates of the block by \(I_k\)
in increasing order. Let \(E_k=\ell_q(I_k)\) and let \(Y_k\) be the same
coordinate space with norm \(\|z\|_{Y_k}=\max\{\|z\|_q,b_k\|z\|_2\}\).
Define the extension by zero \(j_k:E_k\to X_{q,w}\) and the restriction
\(P_k:X_{q,v}\to Y_k\) by
\[
 (j_kx)_i=
 \begin{cases}
 x_i,& i\in I_k,\\
 0,& i\notin I_k,
 \end{cases}
 \qquad P_ky=(y_i)_{i\in I_k}.
\]
Since \(a_k=N_k^{1/q-1/2}\) and
\(\|x\|_2\le N_k^{1/2-1/q}\|x\|_q\) on \(E_k\), these maps satisfy
\[
\begin{aligned}
 \|j_kx\|_{q,w}&=\max\{\|x\|_q,a_k\|x\|_2\}=\|x\|_q,\\
 \|P_ky\|_{Y_k}&=\max\{\|y|_{I_k}\|_q,\|(vy)|_{I_k}\|_2\}
 \le \max\{\|y\|_q,\|vy\|_2\}=\|y\|_{q,v}.
\end{aligned}
\]
For \(x\in E_k\) and \(i\in I_k\), \((P_kDj_kx)_i=d_i x_i=\delta_k x_i\).
Thus \(P_kDj_k=D_{N_k}:E_k\to Y_k\). Setting
\(S_k=P_kRj_k=(P_kB)(Aj_k)\), we obtain
\[
\left\|D_{N_{k}}-S_{k}\right\|=\left\|P_{k}(D-B A) j_{k}\right\| \leq \varepsilon, \quad \gamma_{s}\left(S_{k}\right) \leq\left\|P_{k} B\right\|\left\|A j_{k}\right\| \leq\|B\|\|A\| .
\]

Set $\alpha=1 / s-1 / q>0$. By Proposition~\ref{res:2}, every block satisfies
\[
(1-\varepsilon) \frac{2^{\alpha k}}{k} \leq(1-\varepsilon) \max \left\{1, k^{-1} 2^{\alpha k}\right\} \leq \gamma_{s}\left(S_{k}\right) \leq M \quad\left(k \geq k_{0}\right) .
\]

Both $M$ and $\varepsilon$ are fixed, since every $S_{k}$ is a compression of the same
factorization $R=B A$. To choose a block that violates this bound, put $u=2^{\alpha}-1>0$. For
$k \geq 2$, the binomial formula gives
\[
\frac{2^{\alpha k}}{k}=\frac{(1+u)^{k}}{k} \geq \frac{\binom{k}{2} u^{2}}{k}=\frac{(k-1) u^{2}}{2} .
\]

Choose an integer $k \geq \max \left\{2, k_{0}\right\}$ with $k>1+2 M /\left((1-\varepsilon) u^{2}\right)$. For this block,
\[
M<(1-\varepsilon) \frac{(k-1) u^{2}}{2} \leq(1-\varepsilon) \frac{2^{\alpha k}}{k} \leq \gamma_{s}\left(S_{k}\right) \leq M,
\]
a contradiction. Hence $\|D-R\| \geq 1$ for every $R \in \mathcal{G}_{s}\left(X_{q, w}, X_{q, v}\right)$.
If $R_{n} \in \mathcal{G}_{s}\left(X_{q, w}, X_{q, v}\right)$ and $R_{n} \rightarrow R$ in
operator norm, then $\|D-R\|=\lim _{n \rightarrow \infty}\left\|D-R_{n}\right\| \geq$ 1. This
proves the lower bound in \eqref{eq:7} for the whole closure. The zero operator factors
through $L_{s}$ and satisfies $\|D-0\|=\|D\|=1$, giving equality.
\end{proof}

For example, when $p=4 / 3$ and $r=3 / 2$, the conjugate exponents are $q=4$ and $s=3$.
Among approximants to the $k$-th block with error at most $1/2$, the minimum
product of the factor norms through $L_3$ is
\[
\frac{1}{2} \max \left\{1,2^{k / 12} / k\right\} .
\]

\subsection{Complemented embeddings in \texorpdfstring{$L_{q}$}{Lq} and duality}
We use a complemented embedding \(J_v\colon X_{q,v}\to L_q\) and a bounded
surjection \(Q_w\colon L_q\to X_{q,w}\) to form \(J_v D Q_w\) on \(L_q\).
Its adjoint acts on \(L_p\). Rosenthal's construction \cite{Ros} yields
these maps and their one-sided inverses with bounds independent of the
weights, as required for \eqref{eq:1}.

\begin{proposition}\label{res:5}
For each $2<q<\infty$, there is a constant $C_q>0$ such that for every sequence
$u=(u_j)$ with $0<u_j\leq 1$ there are bounded linear maps
$J_u:X_{q,u}\to L_q$ and $Q_u:L_q\to X_{q,u}$ satisfying
\[
 Q_u J_u=I_{X_{q,u}},\qquad \|Q_u\|\leq 1,\qquad \|J_u\|\leq C_q.
\]
In particular, $J_u Q_u$ is a bounded projection onto $J_u(X_{q,u})$.
\end{proposition}
\begin{proof}
Set $q^{\prime}=q /(q-1)$ and $\lambda_{j}=u_{j}^{2 q /(q-2)}$. Choose independent symmetric
random variables $g_{j}$ on (0, 1), taking the values $\pm \lambda_{j}^{-1 / q}$ with
probability $\lambda_{j} / 2$ each and zero otherwise. They can be realized on the atomless
Lebesgue space by using independent uniform random variables. Since
\(\lvert g_j\rvert=\lambda_j^{-1/q}\) on a set of measure \(\lambda_j\) and vanishes
elsewhere, for \(1\le t<\infty\) we have
\[
 \|g_j\|_t^t=\int_0^1 |g_j(\omega)|^t\,d\omega
 =\lambda_j(\lambda_j^{-1/q})^t,\qquad
 \|g_j\|_t=\lambda_j^{1/t-1/q}.
\]
Taking \(t=q,2,q'\) and substituting \(\lambda_j=u_j^{2q/(q-2)}\) gives
\[
\left\|g_{j}\right\|_{q}=1, \quad\left\|g_{j}\right\|_{2}=\lambda_{j}^{1 / 2-1 / q}=u_{j}, \quad\left\|g_{j}\right\|_{q^{\prime}}=\lambda_{j}^{1 / q^{\prime}-1 / q}=u_{j}^{2} .
\]

For every finitely supported \(x\), Rosenthal's inequality \cite[Theorem 3]{Ros}
applied to the independent mean-zero variables \(x_jg_j\) gives
\[
\left\|\sum_{j} x_{j} g_{j}\right\|_{q} \leq C_{q} \max \left\{\left(\sum_{j}\left|x_{j}\right|^{q}\right)^{1 / q},\left(\sum_{j} u_{j}^{2}\left|x_{j}\right|^{2}\right)^{1 / 2}\right\} .
\]

For complex scalars the same bound follows by applying the real inequality to the real and
imaginary parts, with a change in $C_{q}$. Consequently the synthesis map $J_{u} x=\sum_{j}
x_{j} g_{j}$ extends to $X_{q, u}$ with norm at most $C_{q}$.

For $f \in L_{q}$, define the coefficients
\[
a_{j}(f)=u_{j}^{-2} \int_{0}^{1} f(t) g_{j}(t) d t .
\]

Write \(a=(a_j(f))_j\). To show that \(Q_uf=a\) defines a contraction
into \(X_{q,u}\), we estimate the two terms in its norm,
\(\|(u_ja_j(f))_j\|_{\ell_2}\) and \(\|a\|_{\ell_q}\).
Use the Hilbert inner product $(f, h)_{2}=\int_{0}^{1} f(t) \overline{h(t)} d t$. The
functions $g_{j} / u_{j}$ are orthonormal in $L_{2}$ by independence, centering and the
identity $\left\|g_{j}\right\|_{2}=u_{j}$. Since each $g_{j}$ is real-valued, $u_{j}
a_{j}(f)=\left(f, g_{j} / u_{j}\right)_{2}$ over either scalar field. Bessel's inequality and
$L_{q} \subseteq L_{2}$ now give
\[
\left(\sum_{j} u_{j}^{2}\left|a_{j}(f)\right|^{2}\right)^{1 / 2} \leq\|f\|_{2} \leq\|f\|_{q} .
\]

To estimate the $\ell_{q}$-norm of the coefficients, set $h_{j}=g_{j} / u_{j}^{2}$. These
variables are independent and symmetric, and $\left\|h_{j}\right\|_{q^{\prime}}=1$. For a
finitely supported scalar sequence $b=\left(b_{j}\right)$, the sums below run over its
support. Let \(\sigma_j\) be independent uniform signs. For every fixed
choice of signs, \((\sigma_jh_j)_j\) has the same joint law as \((h_j)_j\).
Averaging over the signs gives
\[
 \Big\|\sum_j b_jh_j\Big\|_{q'}^{q'}
 =\int_0^1\mathbb E_\sigma\Big|\sum_j\sigma_j b_jh_j(t)\Big|^{q'}\,dt .
\]
For fixed \(t\), put \(c_j=b_jh_j(t)\). Since \(0<q'/2<1\), Jensen's
inequality for the concave function \(v\mapsto v^{q'/2}\) gives
\[
 \mathbb E_\sigma\Big|\sum_j\sigma_jc_j\Big|^{q'}
 \le \Big(\mathbb E_\sigma\Big|\sum_j\sigma_jc_j\Big|^2\Big)^{q'/2}
 =\Big(\sum_j|c_j|^2\Big)^{q'/2}.
\]
The equality follows from \(\mathbb E_\sigma(\sigma_i\sigma_j)=0\) for
\(i\ne j\) and \(\sigma_j^2=1\). Subadditivity of \(v\mapsto v^{q'/2}\)
on \([0,\infty)\) gives \((\sum_j|c_j|^2)^{q'/2}\le\sum_j|c_j|^{q'}\).
Integrating in \(t\) and using \(\|h_j\|_{q'}=1\), we obtain
\[
 \Big\|\sum_j b_jh_j\Big\|_{q'}^{q'}
 \le \int_0^1\sum_j |b_jh_j(t)|^{q'}\,dt
 =\sum_j |b_j|^{q'} .
\]

H\"older's inequality therefore gives
\[
\left|\sum_{j} b_{j} a_{j}(f)\right|=\left|\int_{0}^{1} f(t) \sum_{j} b_{j} h_{j}(t) d t\right| \leq\|f\|_{q}\|b\|_{q^{\prime}} .
\]

Taking the supremum over finitely supported $b$ with $\|b\|_{q^{\prime}} \leq 1$ yields
$\|a\|_{\ell_q}\leq\|f\|_q$. The two bounds imply
\[
\left\|a \mathbf{1}_{\{j>n\}}\right\|_{q, u}=\max \left\{\left(\sum_{j>n}\left|a_{j}(f)\right|^{q}\right)^{1 / q},\left(\sum_{j>n} u_{j}^{2}\left|a_{j}(f)\right|^{2}\right)^{1 / 2}\right\} \longrightarrow 0 .
\]

Thus $Q_{u} f=a \in X_{q, u}$ and $\left\|Q_{u} f\right\|_{q, u} \leq\|f\|_{q}$. For finitely supported $x$, orthogonality gives
\[
\left(Q_{u} J_{u} x\right)_{j}=\sum_{i} x_{i} u_{j}^{-2} \int_{0}^{1} g_{i}(t) g_{j}(t) d t=x_{j} .
\]

By density $Q_{u} J_{u}=I_{X_{q, u}}$, so $\|x\|_{q, u} \leq\left\|J_{u} x\right\|_{q}$. Also
$\left(J_{u} Q_{u}\right)^{2}=J_{u} Q_{u}$ and $J_{u} Q_{u}\left(J_{u} x\right)=J_{u} x$,
giving the claimed projection onto $J_{u}\left(X_{q, u}\right)$.
\end{proof}

\begin{proof}[Proof of Theorem~\ref{res:1}]
Fix $1<p<2$ and let $q=p /(p-1)$. Take the weights $w, v$ from \eqref{eq:6} and the
corresponding diagonal operator $D: X_{q, w} \rightarrow X_{q, v}$. With the maps from
Proposition~\ref{res:5}, set
\[
H=J_{v} D Q_{w}: L_{q} \longrightarrow L_{q}, \quad T_{p}=H^{*}: L_{p} \longrightarrow L_{p} .
\]

We identify $L_{q}^{*}$ with $L_{p}$ via the bilinear pairing. By reflexivity, $T_{p}^{*}=H$. Proposition~\ref{res:5} then gives
\begin{equation}
Q_{v} T_{p}^{*} J_{w}=\left(Q_{v} J_{v}\right) D\left(Q_{w} J_{w}\right)=D . \label{eq:8}
\end{equation}

We first show that $T_{p}$ is $\ell_{p}$-strictly singular. Suppose it is bounded below on a
closed subspace $X \cong \ell_{p}$. Its image $T_{p} X$ is then a closed copy of $\ell_{p}$ in
$L_{p}$. By the complemented-subspace property recalled in \cite[p. 1106]{JS}, $T_{p} X$
contains a further copy $Y \cong \ell_{p}$ onto which there is a bounded projection $P: L_{p}
\rightarrow Y$. Choose an isomorphism $W: \ell_{p} \rightarrow Y$. The operators
\[
U=\left(\left.T_{p}\right|_{X}\right)^{-1} W: \ell_{p} \rightarrow L_{p}, \quad V=W^{-1} P: L_{p} \rightarrow \ell_{p}
\]
satisfy $V T_{p} U=W^{-1} P W=I_{\ell_{p}}$, since $W\left(\ell_{p}\right)=Y$ and $\left.P\right|_{Y}=I_{Y}$. Taking adjoints gives
\[
I_{\ell_{q}}=U^{*} T_{p}^{*} V^{*}=\left(U^{*} J_{v}\right) D\left(Q_{w} V^{*}\right) .
\]

Put $M=Q_{w} V^{*}$ and $K=U^{*} J_{v}$. From $K D M=I_{\ell_{q}}$ we obtain
\[
\frac{\|z\|_{q}}{\|K\|\|D\|} \leq\|M z\|_{q, w} \leq\|M\|\|z\|_{q}, \quad\|D(M z)\|_{q, v} \geq \frac{\|z\|_{q}}{\|K\|} \geq \frac{\|M z\|_{q, w}}{\|K\|\|M\|} .
\]

The first estimate identifies $M\left(\ell_{q}\right)$ as a closed copy of $\ell_{q}$. The second bounds $D$ below on that copy, contradicting Lemma~\ref{res:3}.

Now fix $p<r<2$ and write $s=r /(r-1)$, so $2<s<q$. For any $R \in \mathcal{G}_{r}\left(L_{p},
L_{p}\right)$, choose $A: L_{p} \rightarrow L_{r}$ and $B: L_{r} \rightarrow L_{p}$ with $R=B
A$. The compression of $R^{*}$ factors through $L_{s}$, since
\[
Q_{v} R^{*} J_{w}=\left(Q_{v} A^{*}\right)\left(B^{*} J_{w}\right), \quad B^{*} J_{w}: X_{q, w} \rightarrow L_{s}, \quad Q_{v} A^{*}: L_{s} \rightarrow X_{q, v} .
\]

By Lemma~\ref{res:4}, this operator is at distance at least one from $D$. Using \eqref{eq:8} and the norm bounds from Proposition~\ref{res:5}, we obtain
\[
1 \leq\left\|D-Q_{v} R^{*} J_{w}\right\|=\left\|Q_{v}\left(T_{p}^{*}-R^{*}\right) J_{w}\right\| \leq\left\|Q_{v}\right\|\left\|J_{w}\right\|\left\|T_{p}-R\right\| \leq C_{q}\left\|T_{p}-R\right\| .
\]

By continuity, $1 \leq C_{q}\left\|T_{p}-R\right\|$ also holds for every $R \in
\overline{\mathcal{G}_{r}\left(L_{p}, L_{p}\right)}$. The block index in Lemma~\ref{res:4} may
depend on $s=r^{\prime}$, but its distance bound is one for every $2<s<q$. Thus the same
$T_{p}$ and $\eta_{p}=C_{q}^{-1}$ give \eqref{eq:1} for all $p<r<2$.
\end{proof}

\section{Factorization thresholds and the chain of closed ideals}\label{sec:thresholds}
We now vary the block coefficients from Section~\ref{sec:main} to
prescribe the threshold for factorization and determine whether
factorization holds at the threshold.

Fix $2<q<\infty$. Let $I_{k}$ be finite consecutive blocks of lengths
$N_{k}$, put $a_{k}=N_{k}^{1 / q-1 / 2}$, and choose $a_{k} \leq \delta_{k} \leq 1$. On
$I_{k}$ set $w_{j}=a_{k}, d_{j}=\delta_{k}$, and $v_{j}=w_{j} / d_{j}$. As in Subsection 2.2,
the diagonal $D x=d x$ is a norm-one operator from $X_{q, w}$ to $X_{q, v}$.
\begin{theorem}\label{res:6}
For $2<s<q$, let $M_{s}=\sup _{k} \delta_{k} N_{k}^{1 / s-1 / q}$. If $M_{s}<\infty$, then $D$ factors through $L_{s}$. If $M_{s}=\infty$, then
\[
\operatorname{dist}\left(D, \overline{\mathcal{G}_{s}\left(X_{q, w}, X_{q, v}\right)}\right)=1 .
\]

In particular, $D$ factors through $L_s$ if and only if $D\in\overline{\mathcal G_s(X_{q,w},X_{q,v})}$.

\end{theorem}
\begin{proof}
Suppose $M_{s}=\infty$. If $S=B A$ factors through $L_{s}$ and $\|D-S\| \leq \varepsilon<1$, compress to each block as in Lemma~\ref{res:4}. Proposition~\ref{res:2} gives
\[
(1-\varepsilon) \delta_{k} N_{k}^{1 / s-1 / q} \leq \gamma_{s}\left(P_{k} S j_{k}\right) \leq\|B\|\|A\| \quad \text { for all } k \text {. }
\]

Taking the supremum over \(k\) gives
\[
 M_s=\sup_k \delta_k N_k^{1/s-1/q}
 \leq \frac{\|A\|\,\|B\|}{1-\varepsilon}<\infty.
\]
This contradicts the assumption \(M_s=\infty\). Thus every operator factoring through
\(L_s\) is at distance at least one from \(D\). By continuity, the same bound
holds throughout the norm closure of this class.
The zero operator attains this distance since \(\|D\|=1\).

Suppose $M_{s}<\infty$, and choose $0<\theta<1$ with $1 / s=(1-\theta) / q+\theta / 2$. Write $t_{j}=w_{j}^{\theta}$ and $z_{j}=w_{j}^{1-\theta}$. Consider the maps
\[
X_{q, w} \xrightarrow{A} X_{s, z} \xrightarrow{B} X_{q, v}, \quad(A x)_{j}=t_{j} x_{j}, \quad(B y)_{j}=\left(d_{j} / t_{j}\right) y_{j} .
\]

For finitely supported $x$, apply H\"older's inequality to $\left|t_{j}
x_{j}\right|^{s}=\left|x_{j}\right|^{(1-\theta) s}\left|w_{j} x_{j}\right|^{\theta s}$ with
conjugate exponents $q /((1-\theta) s)$ and $2 /(\theta s)$. This gives
\[
\|t x\|_{s} \leq\|x\|_{q}^{1-\theta}\|w x\|_{2}^{\theta}, \quad\|z t x\|_{2}=\|w x\|_{2} .
\]

Thus $A$ extends as a contraction. For $j \in I_{k}, t_{j}=N_{k}^{1 / q-1 / s}, d_{j} / t_{j} \leq M_{s}$, and $v_{j} d_{j} / t_{j}=z_{j}$. Consequently
\[
\|B y\|_{q, v}=\max \left\{\|(d / t) y\|_{q},\|z y\|_{2}\right\} \leq \max \left\{M_{s}\|y\|_{s},\|z y\|_{2}\right\} \leq \max \left\{1, M_{s}\right\}\|y\|_{s, z} .
\]

On finitely supported vectors $B A x=d x$, so density gives $D=B A$. Proposition~\ref{res:5},
applied with exponent $s$, gives maps $J_{z}: X_{s, z} \rightarrow L_{s}$ and $Q_{z}: L_{s}
\rightarrow X_{s, z}$ with $Q_{z} J_{z}=I$. The required factorization is $D=\left(B
Q_{z}\right)\left(J_{z} A\right)$.
\end{proof}

\begin{theorem}\label{res:7}
Fix $1<p<r_{0}<2$ and let $U$ be either $\left(p, r_{0}\right]$ or $\left(p, r_{0}\right)$.
There is an $\ell_{p}$-strictly singular operator $T: L_{p} \rightarrow L_{p}$ such that $T
\in \mathcal{G}_{r}\left(L_{p}, L_{p}\right)$ for every $r \in U$ and
\[
\operatorname{dist}\left(T, \overline{\mathcal{G}_{r}\left(L_{p}, L_{p}\right)}\right) \geq C_{p^{\prime}}^{-1} \quad(r \in(p, 2) \backslash U),
\]
where $C_{p'}$ is the constant given by Proposition~\ref{res:5} with $q=p'$.
\end{theorem}
\begin{proof}
Put $q=p^{\prime}, \alpha_{0}=1 / p-1 / r_{0}$, and $N_{k}=2^{k}$. For $U=\left(p,
r_{0}\right]$, take $\delta_{k}=N_{k}^{-\alpha_{0}}$. For $U=\left(p, r_{0}\right)$, take
$\delta_{k}=k N_{k}^{-\alpha_{0}}$. Since $0<\alpha_{0}<1 / 2-1 / q$, both choices satisfy
$a_{k} \leq \delta_{k} \leq 1$ after discarding finitely many blocks, and $\delta_{k}
\rightarrow 0$. For $s=r^{\prime}$, the quantity $\delta_{k} N_{k}^{1 / s-1 / q}$ equals
$2^{k\left(1 / r_{0}-1 / r\right)}$ in the first case and $k 2^{k\left(1 / r_{0}-1 /
r\right)}$ in the second. In each case its supremum is finite exactly for $r \in U$.

As in the proof of Theorem~\ref{res:1}, set $H=J_{v} D Q_{w}$ and $T=H^{*}$. The proof of
Lemma~\ref{res:3} uses only $d_{j} \rightarrow 0$ and $v_{j} d_{j}=w_{j}$, so $D$ is
$\ell_{q}$-strictly singular. If $T$ preserved a copy of $\ell_{p}$, the factorization of
$I_{\ell_{q}}$ through $D$ in the proof of Theorem~\ref{res:1} would contradict
Lemma~\ref{res:3}. For $r \in U$, write $D=B A$ through $L_{s}$, using Theorem~\ref{res:6}.
Then
\[
T=\left(J_{v} B A Q_{w}\right)^{*}=\left(Q_{w}^{*} A^{*}\right)\left(B^{*} J_{v}^{*}\right) \in \mathcal{G}_{r}\left(L_{p}, L_{p}\right), \quad L_{s}^{*}=L_{r} .
\]

For $p<r<2$ outside $U$, each $R \in \mathcal{G}_{r}\left(L_{p}, L_{p}\right)$ satisfies
$Q_{v} R^{*} J_{w} \in \mathcal{G}_{s}\left(X_{q, w}, X_{q, v}\right)$. By Theorem~\ref{res:6}
and \eqref{eq:8},
\[
1 \leq\left\|D-Q_{v} R^{*} J_{w}\right\|=\left\|Q_{v}\left(T^{*}-R^{*}\right) J_{w}\right\| \leq C_{q}\|T-R\| .
\]

Taking the infimum in $R$ gives the same bound for the closure.
\end{proof}

For fixed \(1<p<2\), put \(J_r=\overline{\mathcal G_r(L_p,L_p)}\) for
\(p<r\le 2\). Each \(J_r\) is a closed two-sided ideal. Indeed,
\(L_r\oplus_r L_r\) is isometric to \(L_r\), so sums of two
\(L_r\)-factorizations again factor through \(L_r\).

\begin{corollary}\label{cor:jumps}
For every \(p<a<2\),
\[
  \overline{\bigcup_{a<r\le 2}J_r}
  \ \subsetneq\ J_a\ \subsetneq\
  \bigcap_{p<r<a}J_r .
\]
\end{corollary}
\begin{proof}
We use the following consequence of Maurey's factorization theorem
\cite[Th\'eor\`eme~2]{Maurey}. If \(Z\) has Rademacher type \(t\) and
\(1<p<r<t\le 2\), then every bounded map \(B:Z\to L_p\) factors
through \(L_r\). Indeed, by Pisier \cite[Proposition~3, p.~III.5]{PisierType},
Rademacher type \(t\) gives stable type \(r\) for \(r<t\). Maurey's
change of density then gives \(B=M_g V\), where \(V:Z\to L_r\) is
bounded and \(M_g:L_r\to L_p\) is multiplication by a function in
\(L_u\), \(1/u=1/p-1/r\). Over the complex field, apply the real
factorization to \(z\mapsto\operatorname{Re}(Bz)\) on the underlying real
space of \(Z\), and write \(\operatorname{Re}(Bz)=gV_0(z)\).
Then \(V(z)=V_0(z)-iV_0(iz)\) is bounded and complex linear, and
\(gV(z)=\operatorname{Re}(Bz)-i\operatorname{Re}(B(iz))=Bz\).
Taking \(Z=L_t\) and
applying this to the second factor gives \(J_t\subseteq J_r\).
Thus both inclusions in the statement hold.

Apply Theorem~\ref{res:7} with \(r_0=a\) and \(U=(p,a]\).
The resulting operator belongs to \(J_a\) and has distance at least
\(C_{p'}^{-1}\) from each \(J_r\) with \(a<r<2\).
The same bound holds for \(J_2\), since \(J_2\subseteq J_r\) for
each such \(r\). The bound is independent of \(r\), so it also holds
for the closure of the union. This proves strictness on the left.
For the right inclusion, take \(U=(p,a)\) in Theorem~\ref{res:7}.
The resulting operator factors through every \(L_r\) with \(p<r<a\)
and has distance at least \(C_{p'}^{-1}\) from \(J_a\).
\end{proof}

\section{Strictly singular operators outside the factorization ideal}\label{sec:strict}

\subsection{Construction of a strictly singular counterexample}
Let \(\mathcal I_p\) be the closed ideal generated by the intermediate
factorizations. Since the \(J_r\) form a chain by Corollary~\ref{cor:jumps},
\[
  \mathcal I_p
    =\overline{\operatorname{span}\!\left(\bigcup_{p<r<2}\mathcal G_r(L_p,L_p)\right)}
    =\overline{\bigcup_{p<r<2}J_r}.
\]
The uniform bound in Theorem~\ref{res:1} separates \(T_p\) from
\(\mathcal I_p\). To see why \(T_p\) preserves a copy of \(\ell_2\), put
$x_{k}=N_{k}^{-1 / q} \mathbf{1}_{I_{k}}$ and take finitely supported scalars
$c=\left(c_{k}\right)$. The blocks are disjoint and $v_{j} d_{j}=w_{j}$, so
\[
\left\|\sum_{k} c_{k} x_{k}\right\|_{q, w}=\max \left\{\|c\|_{q},\|c\|_{2}\right\}=\|c\|_{2}=\max \left\{\left\|\left(\delta_{k} c_{k}\right)_{k}\right\|_{q},\|c\|_{2}\right\}=\left\|D \sum_{k} c_{k} x_{k}\right\|_{q, v} .
\]

Thus $D=Q_{v} T_{p}^{*} J_{w}$ is an isometry on the closed span of the $x_{k}$. Since $J_{w}$
is an isomorphic embedding and $Q_{v}$ is contractive, $T_{p}^{*}$ is bounded below on the
$J_{w}$-image of that span. On $L_{p}$, an operator and its adjoint are simultaneously
$\ell_{2}$-strictly singular \cite[p. 1106]{JS}. Hence $T_{p}$ also preserves a copy of
$\ell_{2}$. For the strictly singular example, we change the norm that
combines the range blocks.

\begin{theorem}\label{res:8}
For each $1<p<2$, there is a strictly singular operator $T:L_p\to L_p$ such that
\[
 \operatorname{dist}(T,\mathcal I_p)\geq C_{p'}^{-1},
\]
where $C_{p'}$ is the constant given by Proposition~\ref{res:5} with $q=p'$.
\end{theorem}
\begin{proof}
Put $q=p^{\prime}$ and use the parameters in \eqref{eq:6}, in particular $N_{k}=2^{k}$ and
$\delta_{k}=1 / k$. Write $E_{k}=\ell_{q}^{N_{k}}$ and let $Y_{k}$ denote the corresponding
range block from Proposition~\ref{res:2}. Replace the global range by $Y=\left(\bigoplus_{k}
Y_{k}\right)_{\ell_{q}}$. The use of an outer direct sum to obtain strict singularity while
retaining the action on individual blocks is related to \cite[Remark 4.3]{JSIdeals}.

The coordinate map $L: X_{q, v} \rightarrow Y, L y=\left(\left.y\right|_{I_{k}}\right)_{k}$, is bounded, since for $y_{k}=\left.y\right|_{I_{k}}$ one has
\[
\|L y\|_{Y}^{q} \leq \sum_{k}\left\|y_{k}\right\|_{q}^{q}+\sum_{k}\left\|\left.(v y)\right|_{I_{k}}\right\|_{2}^{q} \leq\|y\|_{q}^{q}+\|v y\|_{2}^{q} .
\]

Set $\widetilde{D}=L D$. On the vectors above, $\left\|\widetilde{D} \sum_{k} c_{k}
x_{k}\right\|_{Y}=\|c\|_{q}$, since every $D x_{k}$ has block norm one. If $\widetilde{D}$
were bounded below on a copy of $\ell_{q}$, then $D$ would be bounded below on the same copy,
contrary to Lemma~\ref{res:3}.

Every infinite-dimensional closed subspace $Z$ of $Y$ contains a copy of $\ell_{q}$. To see
this, let $m_{0}=0$. Inductively choose a unit vector $y_{n} \in Z$ annihilated by the
projection onto the first $m_{n-1}$ blocks, which has finite rank. Choose $m_{n}>m_{n-1}$ so
that the restriction $u_{n}$ of $y_{n}$ to blocks $m_{n-1}+1, \ldots, m_{n}$ satisfies
$\left\|y_{n}-u_{n}\right\|_{Y}<2^{-n-3}$. The coordinate projections are contractions, so $3
/ 4 \leq\left\|u_{n}\right\|_{Y} \leq 1$. For finitely supported $c$, disjointness gives
\[
\left\|\sum_{n} c_{n} u_{n}\right\|_{Y}^{q}=\sum_{n}\left|c_{n}\right|^{q}\left\|u_{n}\right\|_{Y}^{q}, \quad\left\|\sum_{n} c_{n}\left(y_{n}-u_{n}\right)\right\|_{Y} \leq\|c\|_{\infty} \sum_{n=1}^{\infty} 2^{-n-3} \leq \frac{1}{8}\|c\|_{q} .
\]

Thus $\frac{5}{8}\|c\|_{q} \leq\left\|\sum_{n} c_{n} y_{n}\right\|_{Y} \leq
\frac{9}{8}\|c\|_{q}$, and the closed span of the $y_{n}$ is a copy of $\ell_{q}$ in $Z$. If
$\widetilde{D}$ were bounded below on an infinite-dimensional closed subspace $M$, apply this
construction to $Z=\widetilde{D}(M)$. For the resulting copy $Z_{0} \cong \ell_{q}$, the
subspace $\left(\left.\widetilde{D}\right|_{M}\right)^{-1}\left(Z_{0}\right) \cong \ell_{q}$
contradicts the preceding paragraph. Hence $\widetilde{D}$ is strictly singular.

We next realize $Y$ as a complemented subspace of $L_{q}$. The finite-coordinate version of
Proposition~\ref{res:5} gives $J_{k}: Y_{k} \rightarrow L_{q}$ and $Q_{k}: L_{q} \rightarrow
Y_{k}$ with $\left\|J_{k}\right\| \leq C_{q},\left\|Q_{k}\right\| \leq 1$, and
$Q_{k} J_{k}=I$. Partition (0,1), up to endpoints, into intervals $\left(t_{k}, t_{k}+\lambda_{k}\right)$ of lengths $\lambda_{k}>0$. The rescaling
\[
(\Phi f)_{k}(u)=\lambda_{k}^{1 / q} f\left(t_{k}+\lambda_{k} u\right), \quad 0<u<1, \quad \sum_{k}\left\|(\Phi f)_{k}\right\|_{q}^{q}=\|f\|_{q}^{q}
\]
is onto the $\ell_{q}$-sum of copies of $L_{q}$, by rescaling each component back onto its
interval. With $J_{Y}=\Phi^{-1}\left(\bigoplus_{k} J_{k}\right)$ and
$Q_{Y}=\left(\bigoplus_{k} Q_{k}\right) \Phi$, we have
\[
Q_{Y} J_{Y}=\bigoplus_{k} Q_{k} J_{k}=I_{Y}, \quad\left\|J_{Y}\right\| \leq \sup
_{k}\left\|J_{k}\right\| \leq C_{q}, \quad\left\|Q_{Y}\right\| \leq \sup
_{k}\left\|Q_{k}\right\| \leq 1 .
\]

Put $H=J_{Y} \widetilde{D} Q_{w}$ and $T=H^{*}$. The ideal property of strictly singular
operators gives strict singularity of $H$. On $L_{q}(0,1), 1<q<\infty$, strict singularity is
equivalent to strict singularity of the adjoint \cite[p. 1106]{JS}. Thus $T$ is strictly
singular on $L_{p}$.

Let $j_{k}: E_{k} \rightarrow X_{q, w}$ be the block inclusion and $P_{k}: Y \rightarrow
Y_{k}$ the coordinate projection. Put $U_{k}=J_{w} j_{k}$ and $V_{k}=P_{k} Q_{Y}$. Their norms
are at most $C_{q}$ and one, respectively. For $R \in \mathcal{L}(L_p)$, its adjoint $R^{*}:
L_{q} \rightarrow L_{q}$ gives a block map $V_{k} R^{*} U_{k}: E_{k} \rightarrow Y_{k}$. In
the common coordinates indexed by $I_{k}$, define
\[
\Theta_{k}(R)=\frac{\operatorname{tr}\left(V_{k} R^{*} U_{k}\right)}{N_{k} \delta_{k}}=\frac{1}{N_{k} \delta_{k}} \sum_{j \in I_{k}}\left(V_{k} R^{*} U_{k} e_{j}\right)_{j} .
\]

By the averaging argument of Proposition~\ref{res:2} and the identity
$\|I:E_k\to Y_k\|=\delta_k^{-1}$, we have $\left|\Theta_{k}(R)\right| \leq \left\|V_{k} R^{*}
U_{k}\right\| \leq C_{q}\|R\|$. Since $Q_{Y} T^{*} J_{w}=\widetilde{D}$, the compression of
$T^{*}$ is $V_{k} T^{*} U_{k}=\delta_{k} I$ and $\Theta_{k}(T)=1$. The map $R
\mapsto\left(\Theta_{k}(R)\right)_{k}$ is therefore a bounded linear map from
$\mathcal{L}(L_p)$ to $\ell_{\infty}$, of norm at most $C_{q}$. The preimage
$\mathcal{N}=\left\{R: \Theta_{k}(R) \rightarrow 0\right\}$ of $c_{0}$ is closed and linear.
Since $\Theta_{k}(T)=1$ for every $k$,
\begin{equation}
1=\lim _{k}\left|\Theta_{k}(T-R)\right| \leq C_{q}\|T-R\| \quad(R \in \mathcal{N}) . \label{eq:9}
\end{equation}

If $R=B A$ factors through $L_{r}, p<r<2$, the compression $\left(V_{k}
A^{*}\right)\left(B^{*} U_{k}\right)$ factors through $L_{s}, s=r^{\prime}$, with norm product
at most $C_{q}\|A\|\|B\|$. Applying \eqref{eq:5} and dividing by $\delta_{k}=1 / k$ gives
\[
\left|\Theta_{k}(R)\right| \leq C_{q}\|A\|\|B\| k 2^{-k(1 / s-1 / q)} \longrightarrow 0 .
\]

Thus $\mathcal{I}_{p} \subseteq \mathcal{N}$, and \eqref{eq:9} proves the assertion.
\end{proof}

\subsection{The type threshold and finitely strictly singular perturbations}
By the Maurey factorization used in Corollary~\ref{cor:jumps}, an operator
on \(L_p\) that factors through a space of type \(t>p\) also factors through
\(L_r\) for every \(p<r<t\). Its block traces therefore tend to zero by the
proof of Theorem~\ref{res:8}. The direct estimate below gives the decay rate
in terms of \(t\). For finitely strictly singular operators we bound the
same traces by the arithmetic means of their Bernstein numbers.

A Banach space $Z$ has Rademacher type $t, 1<t \leq 2$, if there is a constant $K_{t}(Z)<\infty$ such that for every finite sequence $z_{1}, \ldots, z_{N}$ in $Z$,
\[
\left(\mathbb{E}_{\sigma}\left\|\sum_{j=1}^{N} \sigma_{j} z_{j}\right\|^{2}\right)^{1 / 2} \leq K_{t}(Z)\left(\sum_{j=1}^{N}\left\|z_{j}\right\|^{t}\right)^{1 / t},
\]
where the $\sigma_j$ are independent uniform signs. Write $\mathcal{H}_{>p}$ for the
operators on $L_{p}$ that factor through a Banach space of type $t$ for at least one $t>p$.
Write $\mathcal{F}$ for the finitely strictly singular operators on $L_{p}$. For a bounded map
$S: E \rightarrow F$ write
\[
b_{m}(S)=\sup _{\substack{H \subseteq E \\ \dim H=m}} \inf _{\substack{x \in H \\\|x\|_{E}=1}}\|S x\|_{F}
\]
whenever $m \leq \operatorname{dim} E$. Finite strict singularity on an infinite-dimensional
domain means $b_{m}(S) \rightarrow 0$. For bounded composable maps $A,S,B$, the Bernstein
numbers satisfy
\[
b_{m}(A S B) \leq\|A\|\|B\| b_{m}(S)
\]
whenever both sides are defined. If $\left.B\right|_{H}$ has a nontrivial kernel on an
$m$-dimensional subspace $H$, the lower bound of $\left.A S B\right|_{H}$ is zero. Otherwise
set $F=B(H)$. For unit $y \in F$, take $x=\left(\left.B\right|_{H}\right)^{-1} y$, so $1
\leq\|B\|\|x\|$. It follows that
\[
\inf _{\substack{z \in H \\\|z\|=1}}\|A S B z\| \leq \frac{\|A S y\|}{\|x\|} \leq\|A\|\|B\|\|S y\| .
\]

Taking the infimum over unit $y \in F$ and then the supremum over $H$ gives the stated inequality, since $\operatorname{dim} F=m$.

\begin{theorem}\label{res:9}
Let $1<p<2$ and $q=p /(p-1)$. The strictly singular operator $T$ from Theorem~\ref{res:8} satisfies
\[
\operatorname{dist}\left(T, \overline{\operatorname{span}\left(\mathcal{H}_{>p} \cup \mathcal{F}\right)}\right) \geq C_{q}^{-1},
\]
where $C_q$ is the constant given by Proposition~\ref{res:5}.
\end{theorem}
\begin{proof}
Retain $U_{k}, V_{k}, \Theta_{k}$ and $\mathcal{N}$ from the proof of Theorem~\ref{res:8}.
Write $N=2^{k}$, so $a_{k}= N^{1 / q-1 / 2}$ and the norm on $Y_{k}$ is $\max
\left\{\|\cdot\|_{q},\left(a_{k} / \delta_{k}\right)\|\cdot\|_{2}\right\}$. By \eqref{eq:9} it
suffices to show $\mathcal{H}_{>p} \cup \mathcal{F} \subseteq \mathcal{N}$.

Suppose first that $S=B A$, with $A: L_{p} \rightarrow Z, B: Z \rightarrow L_{p}$, where $Z$
has type $t>p$. Set $F=B^{*} U_{k}: E_{k} \rightarrow Z^{*}$ and $G=V_{k} A^{*}: Z^{*}
\rightarrow Y_{k}$. Write $\left\langle z, z^{*}\right\rangle=z^{*}(z)$ for $z \in Z$ and
$z^{*} \in Z^{*}$. The coordinate functionals of $Y_{k}$ have norm at most one, hence
\[
|\operatorname{tr}(G F)| \leq\|G\| \sum_{j=1}^{N}\left\|F e_{j}\right\|_{Z^{*}} .
\]

For any $\rho>0$, choose $z_{j} \in B_{Z}$ such that $\operatorname{Re}\left\langle z_{j}, F
e_{j}\right\rangle \geq\left\|F e_{j}\right\|-\rho / N$. Since
$\mathbb{E}_{\sigma}\left(\sigma_{i} \sigma_{j}\right)$ equals one for $i=j$ and zero
otherwise,
\[
\sum_{j=1}^{N}\left\langle z_{j}, F e_{j}\right\rangle=\mathbb{E}_{\sigma}\left\langle\sum_{j=1}^{N} \sigma_{j} z_{j}, F\left(\sum_{j=1}^{N} \sigma_{j} e_{j}\right)\right\rangle .
\]

Write $\sigma=\sum_{j} \sigma_{j} e_{j}$. Since $\|\sigma\|_{q}=N^{1 / q}$, the preceding identity and the type inequality yield
\[
\begin{aligned}
\sum_{j=1}^{N}\left\|F e_{j}\right\|-\rho & \leq \mathbb{E}_{\sigma}\left[\left\|\sum\nolimits_{j} \sigma_{j} z_{j}\right\|\|F \sigma\|\right] \\
& \leq\|F\| N^{1 / q}\left(\mathbb{E}_{\sigma}\left\|\sum\nolimits_{j} \sigma_{j} z_{j}\right\|^{2}\right)^{1 / 2} \\
& \leq\|F\| N^{1 / q} K_{t}(Z)\left(\sum\nolimits_{j}\left\|z_{j}\right\|^{t}\right)^{1 / t} \\
& \leq\|F\| K_{t}(Z) N^{1 / q+1 / t} .
\end{aligned}
\]

Let $\rho \downarrow 0$. Since $1-1 / q=1 / p$, we have proved
\[
\left|\Theta_{k}(B A)\right| \leq C_{q} K_{t}(Z)\|A\|\|B\| k 2^{-k(1 / p-1 / t)} \longrightarrow 0 .
\]

Since \(1/p-1/t>0\) for each such factorization, \(\mathcal H_{>p}\subseteq\mathcal N\).

Next take $S \in \mathcal{F}$. Identify the $N$ coordinates with the sign cube
$\Omega_{k}=\{-1,1\}^{k}$. Put $r_{i}(\omega)=\omega_{i}$ and $\mathbb{E}_{\omega}
h(\omega)=2^{-k} \sum_{\omega \in \Omega_{k}} h(\omega)$. The coordinate norms on $E_{k}$ and
$Y_{k}$ use counting measure. Use matrices $H: \ell_{2}^{k} \rightarrow \ell_{q}^{N}$ and $W:
Y_{k} \rightarrow \ell_{2}^{k}$ defined by
\[
(H u)(\omega)=N^{-1 / q} \sum_{i=1}^{k} u_{i} r_{i}(\omega), \quad(W y)_{i}=N^{1 / q} \mathbb{E}_{\omega}\left[y(\omega) r_{i}(\omega)\right] .
\]

Let $\iota_{k}: E_{k} \rightarrow Y_{k}$ denote the identity on block coordinates in this
calculation. Orthogonality of the $r_{i}$ gives $W \iota_{k} H=I_{\ell_{2}^{k}}$. Khintchine's
inequality gives $\|H\|\leq K_q$, where $K_q$ depends only on $q$. With the uniform
probability measure on the cube, Bessel's inequality gives
\[
\|W y\|_{2}^{2}=N^{2 / q} \sum_{i=1}^{k}\left|\mathbb{E}_{\omega}\left(y r_{i}\right)\right|^{2} \leq N^{2 / q} \mathbb{E}_{\omega}|y|^{2}=a_{k}^{2}\|y\|_{2}^{2} \leq \delta_{k}^{2}\|y\|_{Y_{k}}^{2} .
\]

Let $P$ range over all signed permutation matrices on this block, with $\mathbb{E}_{P}$ denoting their uniform average. Consider
\[
M_{P}=\delta_{k}^{-1} W P^{-1} V_{k} S^{*} U_{k} P H .
\]

Averaging $P^{-1}\left(V_{k} S^{*} U_{k}\right) P$ over signed permutations gives
$\operatorname{tr}\left(V_{k} S^{*} U_{k}\right) \iota_{k} / N$, as in Proposition~\ref{res:2}
. Hence $\mathbb{E}_{P} M_{P}=\Theta_{k}(S) I_{\ell_{2}^{k}}$. Put
\[
L_{P}=\delta_{k}^{-1} W P^{-1} V_{k}, \quad R_{P}=U_{k} P H .
\]

Then $M_{P}=L_{P} S^{*} R_{P}$, with $\left\|L_{P}\right\| \leq 1$ and $\left\|R_{P}\right\|
\leq C_{q} K_{q}$. Reflexivity gives $M_{P}^{*}=R_{P}^{*} S L_{P}^{*}$. Under our bilinear
pairing the matrix of the Banach adjoint $M_{P}^{*}$ is $M_{P}^{\top}$. Over the complex
field, $M_{P}^{\top}=\mathcal{C} M_{P}^{\dagger} \mathcal{C}$, where $\mathcal{C} z=\bar{z}$
and $M_{P}^{\dagger}$ is the Hilbert adjoint. Thus $M_{P}^{*}$ and $M_{P}$ have the same
singular values over either field. The Hilbert-space min-max formula gives
$s_{m}\left(M_{P}\right)=b_{m}\left(M_{P}^{*}\right)$ for the $m$-th singular value in
decreasing order. For $1 \leq m \leq k$,
\[
s_{m}\left(M_{P}\right)=b_{m}\left(R_{P}^{*} S L_{P}^{*}\right) \leq\left\|R_{P}\right\| b_{m}(S)\left\|L_{P}\right\| \leq C_{q} K_{q} b_{m}(S) .
\]

Taking the trace of $\mathbb{E}_{P} M_{P}=\Theta_{k}(S) I_{\ell_{2}^{k}}$ and using
$|\operatorname{tr} M_P|\leq\sum_{m=1}^k s_m(M_P)$, we obtain
\[
\left|\Theta_{k}(S)\right|=\frac{1}{k}\left|\mathbb{E}_{P} \operatorname{tr} M_{P}\right| \leq \frac{1}{k} \mathbb{E}_{P} \sum_{m=1}^{k} s_{m}\left(M_{P}\right) \leq \frac{C_{q} K_{q}}{k} \sum_{m=1}^{k} b_{m}(S) .
\]

Since $b_m(S)\to 0$, the arithmetic means on the right tend to zero. Thus $S\in\mathcal N$.

Since $\mathcal{N}$ is closed and linear, it contains the closed linear span of the two classes. The conclusion follows from \eqref{eq:9}.
\end{proof}

For the strictly singular operator $S=T$ and the identity permutation $P=I, V_{k} T^{*} U_{k}=
\delta_{k} \iota_{k}$ gives $M_{I}=W \iota_{k} H=I_{\ell_{2}^{k}}$. For any $m$, choose a
block with $k \geq m$. The same factorization of $M_{I}^{*}$ through $T$ gives
\begin{equation}
b_{m}(T) \geq\left(C_{q} K_{q}\right)^{-1} \quad(m \geq 1) . \label{eq:10}
\end{equation}

Mathieu and Tradacete \cite[Theorem 6.2]{MT} proved that every strictly singular
operator on $L_p$ factors through $\ell_p$, a space of type $p$.
Theorem~\ref{res:9} gives a strictly singular operator that cannot even be
approximated in norm by operators factoring through spaces of type $t>p$.

Since the block traces in the proof of Theorem~\ref{res:9} vanish in the
limit on the closed span of the two classes, they induce a map on the
quotient. Write \(\mathcal{SS}(L_p)\) for
the closed ideal of strictly singular operators on \(L_p\), and set
\[
  \mathcal J=\overline{\operatorname{span}
      \bigl(\mathcal H_{>p}\cup\mathcal F\bigr)} .
\]
\begin{corollary}\label{cor:quotient}
For every \(1<p<2\), the Banach space
\[
 \mathcal{SS}(L_p)\big/\bigl(\mathcal{SS}(L_p)\cap\mathcal J\bigr)
\]
contains a complemented subspace isomorphic to \(\ell_\infty/c_0\).
\end{corollary}
\begin{proof}
Use the blocks and the strictly singular operator
\(T=(J_Y\widetilde D Q_w)^*\) of Theorem~\ref{res:8}.
For \(b=(b_k)\in\ell_\infty\), let \(M_b\) multiply coordinates in
\(I_k\) by \(b_k\) on \(X_{q,w}\). Then \(\|M_b\|\le\|b\|_\infty\).
Put
\[
 T_b=(J_Y\widetilde D M_b Q_w)^*.
\]
Since \(Q_wJ_w=I\), we have \(T_b=(J_w M_b Q_w)^*T\), so \(T_b\)
is strictly singular. With the bilinear duality pairing used above,
\(b\mapsto T_b\) is linear, and
\(\|T_b\|\le\|J_Y\|\|\widetilde D\|\|Q_w\|\|b\|_\infty\).
For finitely supported \(b\) the operator \(T_b\) has finite rank.
It follows by truncation that \(T_b\) is compact for \(b\in c_0\);
in particular it belongs to \(\mathcal F\subseteq\mathcal J\).
Passing to the quotients therefore gives a bounded map
\[
 E:\ell_\infty/c_0\longrightarrow
   \mathcal{SS}(L_p)\big/\bigl(\mathcal{SS}(L_p)\cap\mathcal J\bigr),
 \qquad E[b]=[T_b].
\]

By the proof of Theorem~\ref{res:9}, the map
\(S\mapsto(\Theta_k(S))_k\) is bounded into \(\ell_\infty\) and
takes \(\mathcal J\) into \(c_0\). Its restriction to
\(\mathcal{SS}(L_p)\) induces a bounded map \(P\) from the latter
quotient to \(\ell_\infty/c_0\). The block compression is
\(V_kT_b^*U_k=b_k\delta_k I\), whence \(\Theta_k(T_b)=b_k\) for
every \(k\). Thus \(PE=I_{\ell_\infty/c_0}\).
In particular \(E\) is an isomorphism onto its closed range and
\(EP\) is a bounded projection onto that range.
\end{proof}

\section{Finite strict singularity and approximation through \texorpdfstring{$L_{r}$}{Lr}}\label{sec:finite}

\subsection{A compactness criterion for the diagonal model}
By \eqref{eq:10}, the Bernstein numbers of the strictly singular counterexample
are bounded away from zero. To study finite strict singularity in the model of
Theorem~\ref{res:8}, we allow arbitrary diagonal coefficients and
coordinate-dependent range weights within each block.

Fix $1<p<2$, put $q=p^{\prime}$ and $1 / u=1 / 2-1 / q$, and use the dyadic blocks $I_{k}$ and
the domain $X=X_{q, w}$ of Subsection 2.2, with $w_{j}=a_{k}=\left|I_{k}\right|^{-1 / u}$ on
$I_{k}$. Choose arbitrary range weights $0<v_{j} \leq 1$. Give $\mathbb{K}^{I_{k}}$ the norm
$\|y\|_{Y_{k}}=\max \left\{\|y\|_{q},\left\|\left(v_{j} y_{j}\right)_{j \in
I_{k}}\right\|_{2}\right\}$, where $\mathbb{K}=\mathbb{R}$ or $\mathbb{C}$, and put
$Y=\left(\bigoplus_{k} Y_{k}\right)_{\ell_{q}}$. Given scalars $d_{j}$, first define $D$ on
finitely supported vectors by
$(D x)_{j}=d_{j} x_{j}$. Whenever $D$ extends boundedly from $X$ to $Y$, set $T=\left(J_{Y} D
Q_{w}\right)^{*}$ on $L_{p}$, using the complemented realizations from Theorem~\ref{res:8}. In
that case $D=Q_{Y} T^{*} J_{w}$.

\begin{theorem}\label{res:10}
Write
\[
M_{k}=\max \left\{\max _{j \in I_{k}}\left|d_{j}\right|,\left(\sum_{j \in I_{k}}\left|v_{j} d_{j}\right|^{u}\right)^{1 / u}\right\} .
\]

The multiplier $D$ extends to a bounded operator from $X$ to $Y$ if and only if $\sup _{k}
M_{k}<\infty$. In this case finite strict singularity of $T$, compactness of $T$, and the
condition $M_{k} \rightarrow 0$ are equivalent. In particular every finitely strictly singular
operator in this family belongs to $\overline{\mathcal{G}_{r}\left(L_{p}, L_{p}\right)}$ for
all $p<r<2$.

\end{theorem}
\begin{proof}
On one domain block the norm is $\|\cdot\|_{q}$. Write $h_{j}=v_{j} d_{j}$ on $I_{k}$. Since
$1 / 2= 1 / q+1 / u$, H\"older's inequality gives $\|h x\|_{2} \leq\|h\|_{u}\|x\|_{q}$. If $h
\neq 0$, equality is attained at $\left|x_{j}\right|^{q}=\left|h_{j}\right|^{u}
/\|h\|_{u}^{u}$. Testing the other norm at a coordinate where $\left|d_{j}\right|$ is maximal
gives
\[
\left\|D_{k}\right\|=\max \left\{\sup _{\|x\|_{q}=1}\|d x\|_{q}, \sup _{\|x\|_{q}=1}\|h x\|_{2}\right\}=\max \left\{\max _{j \in I_{k}}\left|d_{j}\right|,\|h\|_{u}\right\}=M_{k} .
\]

Let $P_{>h}$ retain the range blocks with $k>h$. For finitely supported $x$,
\[
\left\|P_{>h} D x\right\|_{Y}^{q} \leq\left(\sup _{k>h} M_{k}\right)^{q} \sum_{k>h}\left\|\left.x\right|_{I_{k}}\right\|_{q}^{q} \leq\left(\sup _{k>h} M_{k}\right)^{q}\|x\|_{X}^{q} .
\]

Conversely, on each single block the domain norm is $\|\cdot\|_{q}$, so testing there gives
the reverse bound on the operator norm. Taking all blocks proves the boundedness criterion by
density. For bounded $D$, set $D^{(h)}=\left(I-P_{>h}\right) D$ and $T_{h}=\left(J_{Y} D^{(h)}
Q_{w}\right)^{*}$. These maps have finite rank, and
\[
\left\|T-T_{h}\right\| \leq\left\|J_{Y}\right\|\left\|Q_{w}\right\|\left\|D-D^{(h)}\right\| \leq C_{q} \sup _{k>h} M_{k} .
\]

Thus $M_{k} \rightarrow 0$ implies compactness of $T$. If $T$ is compact, dyadic conditional
expectations $E_{n}$ converge uniformly to the identity on
$\overline{T\left(B_{L_{p}}\right)}$, so $\left\|T-E_{n} T\right\| \rightarrow 0$. For
$m>\operatorname{rank}\left(E_{n} T\right)$, every $m$-dimensional subspace meets
$\operatorname{ker}\left(E_{n} T\right)$ nontrivially, giving $b_{m}(T) \leq\left\|T-E_{n}
T\right\|$. Hence $b_{m}(T) \rightarrow 0$. Each finite-rank map factors through some
$\ell_{r}^{d}$, realized as a complemented step-function subspace of $L_{r}$ as in
Proposition~\ref{res:2}. To prove the remaining implication, suppose $M_{k} \nrightarrow 0$.
We show that $b_{m}(T)$ is bounded away from zero in each of the following two cases.

Suppose first that $\left|d_{j_{i}}\right| \geq \eta>0$ for a coordinate $j_{i}$ in each of
infinitely many distinct blocks $I_{k_{i}}$, listed in increasing order. The map $A: \ell_{q}
\rightarrow X, A e_{i}=e_{j_{i}}$, is bounded: the weighted $\ell_{2}$-part is at most
$\left(\sum_{i} a_{k_{i}}^{u}\right)^{1 / u}\|x\|_{q}$, and $a_{k}^{u}=2^{-k}$. The map $B: Y
\rightarrow \ell_{q}$ taking $y$ to $\left(y_{j_{i}} / d_{j_{i}}\right)_{i}$ has norm at most
$1 / \eta$. Substituting $D=Q_{Y} T^{*} J_{w}$ into $B D A=I_{\ell_{q}}$ and taking adjoints
gives
\[
I_{\ell_{p}}=\left(J_{w} A\right)^{*} T\left(B Q_{Y}\right)^{*}, \quad
1=b_{m}\left(I_{\ell_{p}}\right) \leq\left\|J_{w} A\right\|\left\|B Q_{Y}\right\| b_{m}(T)
\quad(m \geq 1),
\]
where $\|J_w A\|$ and $\|B Q_Y\|$ are independent of $m$.
Now suppose $\max _{j \in I_{k}}\left|d_{j}\right| \rightarrow 0$. Since $M_{k}$ does not tend to zero, there are $\eta>0$ and infinitely many $k$ for which
\[
c_{k}=\left(\sum_{j \in I_{k}}\left|v_{j} d_{j}\right|^{u}\right)^{1 / u} \geq \eta .
\]

On each such block define probability weights $\mu_{j}=\left|v_{j} d_{j}\right|^{u} /
c_{k}^{u}$. On these blocks $\max _{j \in I_{k}} \mu_{j} \leq \left(\max _{j \in
I_{k}}\left|d_{j}\right| / \eta\right)^{u} \rightarrow 0$. Fix $m \geq 1$, put $s=2^{m}$, and
choose a block with $\max _{j} \mu_{j} \leq 1 /(4 s)$. In coordinate order, cut after the last
partial sum at most $\ell / s$ for $\ell=1, \ldots, s-1$. Denote
these cut sums by $t_{\ell}$, with $t_{0}=0$ and $t_{s}=1$. Since successive partial sums differ by at most $1 /(4 s)$,
\[
0 \leq \ell / s-t_{\ell} \leq 1 /(4 s), \quad \frac{3}{4 s} \leq t_{\ell}-t_{\ell-1} \leq \frac{5}{4 s} \quad(1 \leq \ell \leq s) .
\]

Index the groups by the $s$ sign vectors $\sigma \in\{-1,1\}^{m}$; their masses $p_{\sigma}$
thus lie between 1/(2s) and $2 / s$. Write $\sigma(j)$ for the group containing $j$.

Define $H: \ell_{2}^{m} \rightarrow \ell_{q}^{I_{k}}$ by
\[
(H z)_{j}=\mu_{j}^{1 / q} \sum_{i=1}^{m} \sigma(j)_{i} z_{i} .
\]

Grouping the terms by sign vector and applying Khintchine's inequality gives
\[
\|H z\|_{q}^{q}=\sum_{\sigma} p_{\sigma}\left|\sum_{i} \sigma_{i} z_{i}\right|^{q} \leq 2 \mathbb{E}_{\sigma}\left|\sum_{i} \sigma_{i} z_{i}\right|^{q} \leq 2 K_{q}^{q}\|z\|_{2}^{q} .
\]

The coordinate map $L: Y_{k} \rightarrow \ell_{2}^{I_{k}}, L y=\left(v_{j} y_{j}\right)_{j}$,
is contractive. Write $D_{k}$ for the block of $D$. Since $1 / u+1 / q=1 / 2$, the lower group
masses and sign orthogonality give
\[
\left\|L D_{k} H z\right\|_{2}^{2}=c_{k}^{2} \sum_{j \in I_{k}} \mu_{j}\left|\sum_{i=1}^{m} \sigma(j)_{i} z_{i}\right|^{2} \geq \frac{1}{2} c_{k}^{2}\|z\|_{2}^{2} .
\]

Put $G=L D_{k} H$, and let $P$ be the orthogonal projection onto $G\left(\ell_{2}^{m}\right)$ in $\ell_{2}^{I_{k}}$. Then
\[
C=G^{-1} P: \ell_{2}^{I_{k}} \rightarrow \ell_{2}^{m}, \quad C G=I_{\ell_{2}^{m}}, \quad\|C\| \leq \sqrt{2} / c_{k} \leq \sqrt{2} / \eta .
\]

Write $j_{k}$ for the inclusion of the domain block and $P_{k}$ for the projection onto the
range block. Since $D_{k}=P_{k} Q_{Y} T^{*} J_{w} j_{k}$, putting $A_{m}=J_{w} j_{k} H$ and
$B_{m}=C L P_{k} Q_{Y}$ gives
\[
I_{\ell_{2}^{m}}=C L D_{k} H=B_{m} T^{*} A_{m}, \quad I_{\ell_{2}^{m}}=A_{m}^{*} T B_{m}^{*} .
\]

The Bernstein number inequality now yields
\[
1 \leq\left\|A_{m}\right\|\left\|B_{m}\right\| b_{m}(T) \leq \frac{2^{1 / 2+1 / q} K_{q} C_{q}}{\eta} b_{m}(T) \quad(m \geq 1) .
\]

Thus $T$ cannot be finitely strictly singular.
\end{proof}

\subsection{Martingale-difference columns and subsequence approximation}
For operators with martingale-difference columns, truncation gives approximation
through \(L_2\) even for some noncompact operators.

Fix $1<p<2$ and a countable measurable partition $\left(C_{i}\right)_{i \geq 1}$ of $(0,1)$
such that each set $C_{i}$ has positive Lebesgue measure $\left|C_{i}\right|$. Define $J:
\ell_{p} \rightarrow L_{p}$ and $Q: L_{p} \rightarrow \ell_{p}$ by
\[
J e_{i}=\left|C_{i}\right|^{-1 / p} \mathbf{1}_{C_{i}}, \quad(Q f)_{i}=\left|C_{i}\right|^{-1 / p^{\prime}} \int_{C_{i}} f(t) d t
\]

Then $J$ is an isometry, $Q J=I_{\ell_{p}}$, and $J Q$ is conditional expectation onto the
$\sigma$-algebra generated by the partition. In particular, $Q$ is a contraction. For a
bounded operator $U: \ell_{p} \rightarrow L_{p}$, the operator $T=U Q$ then has columns
$f_{i}=T J e_{i}=U e_{i}$. Suppose there is a filtration $\left(\mathcal{A}_{i}\right)_{i \geq
0}$ for which each $f_{i}$ is $\mathcal{A}_{i}$-measurable and $\mathbb{E}\left(f_{i} \mid
\mathcal{A}_{i-1}\right)=0$ almost everywhere on (0, 1). The latter condition means that
$\int_{E} f_{i}(t) d t=0$ for every $E \in \mathcal{A}_{i-1}$. Put $B=\sup _{i \geq
1}\left\|f_{i}\right\|_{p}$ and $h(K)=\sup _{i \geq 1}\left\|f_{i}
\mathbf{1}_{\left\{\left|f_{i}\right|>K\right\}}\right\|_{p}$ for $K>0$.
\begin{theorem}\label{res:11}
Let $1<p<2$, let $U: \ell_{p} \rightarrow L_{p}$ be bounded, and put $T=U Q$ with $Q$ as
above. Suppose that $\left(U e_{i}\right)_{i \geq 1}$ is a martingale difference sequence for
a common filtration. For each $p<r \leq 2$ the following conditions are equivalent.
\begin{itemize}
  \item[(1)] $T$ is $\ell_{p}$-strictly singular.
  \item[(2)] $h(K) \rightarrow 0$ as $K \rightarrow \infty$.
  \item[(3)] $T$ is finitely strictly singular.
  \item[(4)] $T \in \overline{\mathcal{G}_{r}\left(L_{p}, L_{p}\right)}$.
\end{itemize}

More precisely, there is a constant $C_{p}$ depending only on $p$ and, for each $K>0$, an
operator $S_{K}: L_{p} \rightarrow L_{2}$ such that, with $\iota_{2, p}: L_{2} \hookrightarrow
L_{p}$ the canonical inclusion,
\[
\left\|T-\iota_{2, p} S_{K}\right\|_{L_{p} \rightarrow L_{p}} \leq C_{p} h(K), \quad\left\|S_{K}\right\|_{L_{p} \rightarrow L_{2}} \leq K^{1-p / 2} B^{p / 2} .
\]

\end{theorem}
\begin{proof}
By Mathieu and Tradacete \cite[Lemma 6.3]{MT}, condition (1) implies that
$T\left(B_{L_{p}}\right)$ is uniformly $p$-integrable, that is,
\(\sup_{\|x\|_p\le1}\|\mathbf 1_E Tx\|_p\to0\) as \(|E|\to0\).
Each $f_{i}=T J e_{i}$ belongs to this
set. The sets $E_{i, K}=\left\{\left|f_{i}\right|>K\right\}$ have measure at most $B^{p} /
K^{p}$, uniformly in $i$. Uniform $p$-integrability on these sets gives $h(K) \rightarrow 0$,
proving (2).

Let \(E_{i-1}\) denote conditional expectation with respect to \(\mathcal A_{i-1}\).
Set \(t_i=f_i\mathbf 1_{\{|f_i|\leq K\}}\) and \(g_i=t_i-E_{i-1}t_i\). Then
$g_{i}$ and $f_{i}-g_{i}=\left(I-E_{i-1}\right)\left(f_{i}-t_{i}\right)$ are martingale
differences. We have $\left\|f_{i}-g_{i}\right\|_{p} \leq 2 h(K)$. For scalar martingale
differences $d_{i}$ and $1<p<2$ the von Bahr--Esseen inequality gives $\left\|\sum_{i}
d_{i}\right\|_{p} \leq c_{p}\left(\sum_{i}\left\|d_{i}\right\|_{p}^{p}\right)^{1 / p}$ for
finite sums \cite{Pinelis}. Applying the real estimate to the real and imaginary parts gives
the complex version, after increasing $c_{p}$ if necessary. It follows that $\left\|\sum_{i}
x_{i}\left(f_{i}-g_{i}\right)\right\|_{p} \leq 2 c_{p} h(K)\|x\|_{\ell_{p}}$ for all finitely
supported coefficient sequences $x$.

For $j<i$, the function $g_{j}$ is $\mathcal{A}_{i-1}$-measurable and $\mathbb{E}\left(g_{i}
\mid \mathcal{A}_{i-1}\right)=0$. Hence $\left(g_{i}, g_{j}\right)_{2}= \int_{0}^{1} g_{i}(t)
\overline{g_{j}(t)} d t=0$. Since $I-E_{i-1}$ is an orthogonal projection on $L_{2}$,
\[
\left\|g_{i}\right\|_{2}^{2} \leq\left\|t_{i}\right\|_{2}^{2}=\int_{\left\{\left|f_{i}\right| \leq K\right\}}\left|f_{i}(t)\right|^{2} d t \leq K^{2-p}\left\|f_{i}\right\|_{p}^{p} \leq K^{2-p} B^{p} .
\]

Orthogonality gives, for finitely supported $x$,
\[
\left\|\sum_{i} x_{i} g_{i}\right\|_{2}^{2}=\sum_{i}\left|x_{i}\right|^{2}\left\|g_{i}\right\|_{2}^{2} \leq K^{2-p} B^{p}\|x\|_{\ell_{2}}^{2} .
\]

The synthesis $V_{K} x=\sum_{i} x_{i} g_{i}$ therefore extends to $\ell_{2}$. Let $I_{p, 2}:
\ell_{p} \hookrightarrow \ell_{2}$ be the coordinate inclusion and set $S_{K}=V_{K} I_{p, 2}
Q$. Since $Q$ and $I_{p, 2}$ are contractions, $\left\|S_{K}\right\| \leq K^{1-p / 2} B^{p /
2}$. For $f \in L_{p}$, the preceding error estimate extends by density to $Q f \in \ell_{p}$,
giving
\[
\left\|\left(T-\iota_{2, p} S_{K}\right) f\right\|_{p}=\left\|\sum_{i}(Q f)_{i}\left(f_{i}-g_{i}\right)\right\|_{p} \leq 2 c_{p} h(K)\|Q f\|_{\ell_{p}} \leq 2 c_{p} h(K)\|f\|_{p} .
\]

For $p<r \leq 2$, the inclusions on the probability interval give $\iota_{2, p} S_{K}=\iota_{r, p}\left(\iota_{2, r} S_{K}\right)$. Thus (2) implies (4).

To deduce (3) from (2), observe that $b_{m}\left(I_{p, 2}\right) \leq m^{1 / 2-1 / p}$.
Indeed, on an $m$-dimensional subspace of $\ell_{p}$ take an extreme point $z$ of its closed
sup-norm unit ball, which is compact in finite dimension. Put
$F=\left\{j:\left|z_{j}\right|=1\right\}$. Since $z \in c_{0}, F$ is finite and $\sup _{j
\notin F}\left|z_{j}\right|<1$. If $|F|<m$, the subspace contains a nonzero $v$ vanishing on
$F$. For sufficiently small $t>0$, both $z+t v$ and $z-t v$ then belong to the ball,
contradicting extremality. Thus $|F| \geq m$ and $\|z\|_{2}^{2} \leq\|z\|_{p}^{p}$ and
$\|z\|_{p} \geq m^{1 / p}$. It follows that $\|z\|_{2} /\|z\|_{p} \leq\|z\|_{p}^{p / 2-1} \leq
m^{1 / 2-1 / p}$. The ideal inequality for Bernstein numbers and $b_{m}(T)
\leq\left\|T-\iota_{2, p} S_{K}\right\|_{L_{p} \rightarrow L_{p}}+b_{m}\left(\iota_{2, p}
S_{K}\right)$ now give
\[
b_{m}(T) \leq C_{p} h(K)+K^{1-p / 2} B^{p / 2} m^{1 / 2-1 / p} .
\]

Let $m \rightarrow \infty$ and then $K \rightarrow \infty$. Condition (3) implies (1) by
definition. Finally, (4) implies (1) by the type argument in the introduction, which applies
also at $r=2$.
\end{proof}

For a noncompact example, take independent $f_{i}$ with values 2, -1 of probabilities 1/3,
2/3, respectively. Since $\mathbb{E} f_{i}=0$, they are martingale differences for their
natural filtration. Independence and $\mathbb{E}\left|f_{i}\right|^{2}=2$ give, for finitely
supported $x$,
\[
\left\|\sum_{i} x_{i} f_{i}\right\|_{p} \leq\left\|\sum_{i} x_{i} f_{i}\right\|_{2}=\sqrt{2}\|x\|_{2} \leq \sqrt{2}\|x\|_{p} .
\]

Thus $U e_{i}=f_{i}$ defines a bounded $U: \ell_{p} \rightarrow L_{p}$, and $T=U Q$ satisfies
Theorem~\ref{res:11} with $h(2)=0$. Yet $\left\|T J e_{i}-T J e_{j}\right\|_{p}=3(4 / 9)^{1 /
p}$ for $i \neq j$, which precludes compactness.

\begin{corollary}\label{res:12}
Fix $1<p<2$. Let $T: L_{p} \rightarrow L_{p}$ be $\ell_{p}$-strictly singular, and let $A:
\ell_{p} \rightarrow L_{p}$ be bounded. There are indices $i_{1}<i_{2}<\cdots$ such that, for
$W: \ell_{p} \rightarrow \ell_{p}$ defined by $W e_{j}=e_{i_{j}}$, the operator $TAW$ belongs to
$\overline{\mathcal{G}_{r}\left(\ell_{p}, L_{p}\right)}$ for every $p<r \leq 2$.

\end{corollary}
\begin{proof}
Since $e_{i} \rightharpoonup 0$ in $\ell_{p}$, the sequence $f_{i}=T A e_{i}$ is weakly null.
It is uniformly $p$-integrable because it lies in $\|A\| T\left(B_{L_{p}}\right)$, by
\cite[Lemma 6.3]{MT}. For each integer $s \geq 0$, let $P_{s}$ be conditional expectation onto
the $\sigma$-algebra generated by the dyadic intervals of length $2^{-s}$. Thus $P_{s} f$ has
the constant value $2^{s} \int_{I} f(t) d t$ on each such interval $I$, up to null sets. Set
$s_{0}=i_{0}=0$. Inductively choose $i_{j}>i_{j-1}$ with $\left\|P_{s_{j-1}}
f_{i_{j}}\right\|_{p}<2^{-j}$, using the finite rank of $P_{s_{j-1}}$. Then choose
$s_{j}>s_{j-1}$ with $\left\|\left(I-P_{s_{j}}\right) f_{i_{j}}\right\|_{p}<2^{-j}$, using
strong convergence of dyadic expectations. Let $d_{j}=\left(P_{s_{j}}-P_{s_{j-1}}\right)
f_{i_{j}}$. The tower property gives $P_{s_{j-1}} d_{j}=0$, and the errors satisfy
\[
\eta_{j}=f_{i_{j}}-d_{j}=\left(I-P_{s_{j}}\right) f_{i_{j}}+P_{s_{j-1}} f_{i_{j}}, \quad\left\|\eta_{j}\right\|_{p}<2^{1-j} .
\]

For every measurable set $G$ and every integer $a \geq 1$,
\[
\sup _{j>a}\left\|\mathbf{1}_{G} d_{j}\right\|_{p} \leq \sup _{i}\left\|\mathbf{1}_{G} f_{i}\right\|_{p}+2^{-a} .
\]

First make $2^{-a}$ small, then use uniform $p$-integrability of the $f_{i}$ and absolute
continuity for $d_{1}, \ldots, d_{a}$. This proves uniform $p$-integrability of ( $d_{j}$ ).
Put $E_{n} x=\sum_{j=1}^{n} x_{j} \eta_{j}$. By H\"older's inequality these finite-rank maps
converge in norm to $E: \ell_{p} \rightarrow L_{p}$, with
\[
\left\|E-E_{n}\right\| \leq\left(\sum_{j>n}\left\|\eta_{j}\right\|_{p}^{p^{\prime}}\right)^{1 / p^{\prime}} \leq\left(\sum_{j>n} 2^{(1-j) p^{\prime}}\right)^{1 / p^{\prime}} \longrightarrow 0 .
\]

Put $U_{d}=T A W-E: \ell_{p} \rightarrow L_{p}$, so $U_{d} e_{j}=d_{j}$. The operator
$\widehat{T}=U_{d} Q: L_{p} \rightarrow L_{p}$ has columns $\widehat{T} J e_{j}=d_{j}$, since
$Q J=I_{\ell_{p}}$. Write $h_{d}(K)=\sup _{j}\left\|d_{j}
\mathbf{1}_{\left\{\left|d_{j}\right|>K\right\}}\right\|_{p}$. The measures
$\left|\left\{\left|d_{j}\right|>K\right\}\right|$ are bounded by
$K^{-p}\left\|U_{d}\right\|^{p}$, uniformly in $j$, so uniform $p$-integrability gives
$h_{d}(K) \rightarrow 0$. Apply Theorem~\ref{res:11} to $\widehat{T}$ and set $V_{K}=S_{K} J$.
Since $J$ is an isometry,
\[
\left\|T A W-E-\iota_{2, p} V_{K}\right\|=\left\|\left(\widehat{T}-\iota_{2, p} S_{K}\right) J\right\| \leq C_{p} h_{d}(K) \longrightarrow 0 .
\]

It follows that
\[
\left\|T A W-\left(E_{n}+\iota_{2, p} V_{K}\right)\right\| \leq\left\|E-E_{n}\right\|+\left\|T A W-E-\iota_{2, p} V_{K}\right\| \longrightarrow 0 \quad(n, K \rightarrow \infty) .
\]

For every $p<r \leq 2$ the two summands factor through $L_{r}$. If they are $B_{1} A_{1}$ and
$B_{2} A_{2}$, their sum factors as $x \mapsto\left(A_{1} x, A_{2} x\right) \mapsto B_{1}
A_{1} x+B_{2} A_{2} x$ through $L_{r} \oplus_{r} L_{r} \cong L_{r}(0,1)$, using the interval
rescaling in Theorem~\ref{res:8}. The chosen subsequence is independent of $r$.
\end{proof}

\begin{samepage}
These results lead to the following question.

\begin{problem}\label{prob:2}
Let $1<p<r<2$, and suppose that $T: L_{p}(0,1) \rightarrow L_{p}(0,1)$ is finitely strictly
singular. For every $\varepsilon>0$, do there exist bounded operators $A: L_{p}(0,1)
\rightarrow L_{r}(0,1)$ and $B: L_{r}(0,1) \rightarrow L_{p}(0,1)$ such that $\|T-B
A\|<\varepsilon$?
\end{problem}
\end{samepage}

\section*{AI Declaration}
During the preparation of this work, the authors used GPT-6.2 to assist with
language editing and the refinement of mathematical arguments. The authors
take full responsibility for the content of this paper.

\end{document}